\documentclass[preprint,12pt]{elsarticle}
\journal{Journal of Differential Equations}
\usepackage{graphicx} 
\usepackage{amsmath,amstext,amssymb,bm}
\usepackage{leftidx}
\usepackage{listings}

\usepackage{xcolor} 
\usepackage{soul} 
\usepackage{tikz}
\usetikzlibrary{shapes,arrows}
\usepackage{mathrsfs}
\usepackage{epstopdf}
\usepackage{color}
\usepackage{multirow}
\usepackage{tabularx}
\usepackage[shortlabels]{enumitem}

\usepackage{amsthm}

\numberwithin{equation}{section}
\newtheorem{definition}{Definition}[section]

\newtheorem{remark}{Remark}[section]
\newtheorem{theorem}{Theorem}[section]

\newtheorem{lemma}{{Lemma}}[section]

\allowdisplaybreaks[4]

\renewcommand{\div}{\operatorname{div}}
\newcommand{\curl}{\operatorname{curl}}

\newcommand{\eps}{\epsilon}

\definecolor{codegreen}{rgb}{0,0.6,0}
\definecolor{codegray}{rgb}{0.5,0.5,0.5}
\definecolor{codepurple}{rgb}{0.58,0,0.82}
\definecolor{backcolour}{rgb}{0.95,0.95,0.92}

\lstdefinestyle{mystyle}{
	backgroundcolor=\color{backcolour},   
	commentstyle=\color{codegreen},
	keywordstyle=\color{magenta},
	numberstyle=\tiny\color{codegray},
	stringstyle=\color{codepurple},
	basicstyle=\ttfamily\footnotesize,
	breakatwhitespace=false,         
	breaklines=true,                 
	captionpos=b,                    
	keepspaces=true,                 
	numbers=left,                    
	numbersep=5pt,                  
	showspaces=false,                
	showstringspaces=false,
	showtabs=false,                  
	tabsize=2
}

\newcommand{\Integer}{\mbox{${\bf Z}$}}
\newcommand{\Real}{\mbox{${\bf R}$}}

\newcommand{\Grad}[1]{\nabla #1}

\newcommand{\Div}[1]{\div\left[#1\right]}

\newcommand{\Curl}[1]{\curl\left[#1\right]}
\newcommand{\Laplacian}[1]{\Delta #1}

\newcommand{\Norm}[2]{\left\|#1\right\|_{#2}}

\newcommand{\SupNorm}[1]{\left|#1\right|_{L^{\infty}}}
\newcommand{\HolderNorm}[2]{\left|#1\right|_{C^{#2}}}
\newcommand{\SobNorm}[2]{\left\|#1\right\|_{H^{#2}}}

\newcommand{\Abs}[1]{\left|#1\right|}

\newcommand{\ImagPart}[1]{\text{Im\{}#1\text{\}}}

\newcommand{\dftl}[1]{\; d#1}
\newcommand{\dV}{\dftl{V}}
\newcommand{\dS}{\dftl{S}}

\newcommand{\cB}{\mathcal{B}}

\newcommand{\cL}{\mathcal{L}}

\newcommand{\cP}{\mathcal{P}}

\newcommand{\px}{\partial_x}
\newcommand{\py}{\partial_y}
\newcommand{\pz}{\partial_z}

\newcommand{\sump}{\sum_{p=-\infty}^{\infty}}
\newcommand{\sumq}{\sum_{q=-\infty}^{\infty}}

\newcommand{\be}{\begin{equation}}
	\newcommand{\ee}{\end{equation}}
\newcommand{\bes}{\begin{equation*}}
	\newcommand{\ees}{\end{equation*}}
\newcommand{\bse}{\begin{subequations}}
	\newcommand{\ese}{\end{subequations}}

\newcommand{\void}[1]{}

\newcommand{\epsu}{\epsilon^{(u)}}
\newcommand{\epsv}{\epsilon^{(v)}}
\newcommand{\epsw}{\epsilon^{(w)}}
\newcommand{\epsm}{\epsilon^{(m)}}

\newcommand{\gammau}{\gamma^{(u)}}

\newcommand{\gammaw}{\gamma^{(w)}}

\newcommand{\gammam}{\gamma^{(m)}}

\newcommand{\ku}{k^{(u)}}

\newcommand{\bepsilon}{\bar{\epsilon}}

\newcommand{\tepsilon}{\mathcal{E}}

\newcommand{\Es}{E^{\text{scat}}}
\newcommand{\Ei}{E^{\text{inc}}}
\newcommand{\Hs}{H^{\text{scat}}}
\newcommand{\Hi}{H^{\text{inc}}}

\newcommand{\tU}{\tilde{U}}
\newcommand{\tW}{\tilde{W}}

\newcommand{\Esx}{E^{\text{scat},x}}
\newcommand{\Esy}{E^{\text{scat},y}}
\newcommand{\Esz}{E^{\text{scat},z}}

\newcommand{\Hcurl}{H(\curl)}
\newcommand{\Hdiv}{H(\div)}
\newcommand{\Hmhcurl}{H^{-1/2}(\curl)}
\newcommand{\Hmhdiv}{H^{-1/2}(\div)}
\newcommand{\HH}{\mathbb{H}}
\newcommand{\HHp}{\mathbb{H}^{\perp}}
\newcommand{\Hzone}{H^1_0}

\newcommand{\GradGu}[1]{\nabla_{\Gamma_u} #1}
\newcommand{\GradGw}[1]{\nabla_{\Gamma_w} #1}
\newcommand{\GradGm}[1]{\nabla_{\Gamma_m} #1}
\newcommand{\DivGu}[1]{\div_{\Gamma_u}\left[#1\right]}
\newcommand{\DivGw}[1]{\div_{\Gamma_w}\left[#1\right]}
\newcommand{\DivGm}[1]{\div_{\Gamma_m}\left[#1\right]}

\newcommand{\InnerOmega}[2]{\left( #1, #2 \right)_{\Omega}}
\newcommand{\InnerGammau}[2]{\langle #1, #2 \rangle_{\Gamma_u}}
\newcommand{\InnerGammaw}[2]{\langle #1, #2 \rangle_{\Gamma_w}}
\newcommand{\InnerGammam}[2]{\langle #1, #2 \rangle_{\Gamma_m}}

\newcommand{\HcurlNorm}[1]{\left\| #1 \right\|_{H({\curl})}}
\newcommand{\HdivNorm}[1]{\left\| #1 \right\|_{H({\div})}}
\newcommand{\HmhdivNorm}[1]{\left\| #1 \right\|_{H^{-1/2}({\div})}}
\newcommand{\HmhcurlNorm}[1]{\left\| #1 \right\|_{H^{-1/2}({\curl})}}

\newcommand{\XNorm}[1]{\Norm{#1}{X}}

\begin{document}

\begin{frontmatter}

\title{A High--Order Perturbation of Envelopes (HOPE) Method
  for Vector Electromagnetic Scattering by Periodic 
  Inhomogeneous Media: Analytic Continuation\tnoteref{tnr1}}

\tnotetext[tnr1]{D.P.N. gratefully acknowledges support from the 
  National Science Foundation through Grant 
  No.~DMS--2111283.}

\author[1]{David P.\ Nicholls\corref{cor1}}
\ead{davidn@uic.edu}

\author[1]{Liet Vo}
\ead{lietvo@uic.edu}

\affiliation[1]{
  organization = {Department of Mathematics, Statistics and Computer Science, 
  University of Illinois at Chicago},
  addressline = {851 South Morgan Street},
  city = {Chicago},
  state = {Illinois},
  postcode = {60607},
  country = {U.S.A.}
  }

\cortext[cor1]{Corresponding author}
	
	\begin{abstract}
		Electromagnetic waves interacting with three--dimensional
		periodic structures occur in many applications of great
		scientific and engineering interest. These three dimensional
		interactions are extremely complicated and subtle, so it is
		unsurprising that practitioners find their rapid, robust, and
		accurate numerical simulation to be of paramount interest.
		Among the wide array of possible numerical approaches, the
		High--Order Spectral algorithms are often preferred due to
		their surpassing fidelity with a moderate number of unknowns,
		and here we describe an algorithm that fits into this class.
		In addition, we take a perturbative approach to the problem
		which views the deviation of the permittivity from a
		reference value as the deformation and we conduct a regular
		perturbation theory. This work concludes a line of research
		on these methods which began with two-dimensional problems
		governed by the Helmholtz equation and moved to small
		perturbations in the fully three-dimensional vector Maxwell
		equations. We now extend these latter results to large
		(real) perturbations constituting a rigorous analytic
		continuation.
	\end{abstract}
	
	\begin{keyword}
		Linear wave scattering,
		Maxwell equations, 
		inhomogeneous media, 
		layered media, 
		High--Order Spectral methods, 
		High--Order Perturbation of Envelopes methods,
		analytic continuation.
  \MSC 65N35 \sep 78M22 \sep 78A45 \sep 35J25 \sep 35Q60 \sep 35Q86
	\end{keyword}
	

\end{frontmatter}
	
	%
	%
	
	\section{Introduction}
	\label{Sec:Intro}
	
	Electromagnetic waves interacting with three--dimensional
	periodic structures occur in many applications of great
	scientific and engineering interest. Many examples can be
	given from fields as different as 
	surface-enhanced spectroscopy \cite{Moskovits85},
	extraordinary optical transmission \cite{ELGTW98},
	cancer therapy \cite{ESHES06}, and 
	surface plasmon resonance (SPR) biosensing
	\cite{Homola08,JJJLWO13,LJJOO12,NichollsReitichJohnsonOh14}.
	
	Because of their crucial role in these linear scattering applications,
	all of the classical numerical algorithms for the simulation of solutions
	to the governing partial differential equations have been brought to
	bear upon this problem. Among these are 
	the Finite Difference \cite{Strikwerda04,LeVeque07}, 
	Finite Element \cite{Johnson87,Ihlenburg98}, 
	Discontinuous Galerkin \cite{HesthavenWarburton08}, 
	Spectral Element \cite{DevilleFischerMund02}, and Spectral
	\cite{GottliebOrszag77,ShenTang06,ShenTangWang11} methods.
	While these are compelling choices, due to their \textit{volumetric}
	character they require a large number of unknowns
	($N = N_x N_y N_z$ for a three dimensional simulation)
	and require the inversion of large, non--symmetric positive definite
	matrices (of dimension $N \times N$). We point the interested reader
	to \cite{ErnstGander12,MoiolaSpence14} for recent developments.
	
	Focusing on the particular example of SPR sensors
	\cite{Homola08,JJJLWO13,LJJOO12,NichollsReitichJohnsonOh14}
	which is the focus of this work, their utility and ubiquity
	follows from two key properties of an SPR, namely its extremely strong and sensitive response. Quantitatively, over
	the range of tens of nanometers in incident wavelength, the reflected
	energy can reduce from almost 100~\% by a factor of 10 or even 100
	before ascending back to almost 100~\%. Clearly, to approximate such
	a structure with the required fidelity, the numerical algorithm should
	produce surpassingly accurate results in a fast and robust manner.
	For this reason, we will focus upon High--Order Spectral (HOS) methods
	\cite{GottliebOrszag77,ShenTang06,ShenTangWang11} which have exactly
	these features.
	
	In regard to the classical approaches listed above, one standard 
	method for generating SPRs is via homogeneous layers of material,
 and it is clearly wasteful to discretize the bulk of each layer. As a result, most prominent solvers feature interfacial unknowns with the
	understanding that information \textit{inside} a layer can be computed
	from appropriate integral formulas. Boundary element (BEM) 
	\cite{SauterSchwab11} and boundary integral (BIM) 
	\cite{ColtonKress13,Kress14} methods are two popular approaches
	and can produce highly accurate solutions in a fraction of the time
	of their volumetric competitors.
	
	In previous work \cite{Nicholls19b,NichollsVo23} the authors investigated
	a new algorithm which has much in common with these HOS algorithms, 
    but was inspired
	by the ``High--Order Perturbation of Surfaces'' (HOPS) methods 
	\cite{NichollsReitich03a,NichollsReitich03b} which have proven to be so
	useful for layered media. A HOPS scheme views the layer interfaces as
	perturbations of flat ones and then	makes recursive corrections to the
	scattering returns from this exactly solvable configuration \cite{Yeh05}.
	However, our new ``High--Order Perturbation of Envelopes'' (HOPE) schemes
	consider more a general permittivity function, $\epsilon(x,y,z)$, which
	does \textit{not} necessarily have layered structure. Our approach
	follows the lead of	Feng, Lin, and Lorton 
	\cite{FengLinLorton15,FengLinLorton16} who adopted a perturbative 
	philosophy by studying the permittivity as a perturbation of a trivial one,
	e.g.,
	\bes
	\epsilon(x,y,z) = \bepsilon (1 - \delta \tepsilon(x,y,z)),
	\quad
	\bepsilon \in \Real,
	\quad
	\tepsilon(x+d_x,y+d_y,z) = \tepsilon(x,y,z),
	\ees
	where $\tepsilon$ is a permittivity ``envelope.''
	In \cite{Nicholls19b} we focused upon
	the two--dimensional scalar problems of electromagnetic radiation in
	Transverse Electric (TE) or Transverse Magnetic (TM) polarization.
	Building upon this work, in \cite{NichollsVo23} we extended our results
	to the three-dimensional vector electromagnetic case governed by
	the full Maxwell equations. These new methods have computational
	advantages over	volumetric solvers in some configurations (e.g., where 
	the support	of $\tepsilon$ is small or where the set on which 
	$\tepsilon$	significantly changes is small). In particular, we considered an approximate indicator function which modeled the absence/presence 
	of a material.
	
	There were several contributions of \cite{Nicholls19b} including a 
	new, and far--reaching, rigorous analysis. In more detail, we proved
	not only that the domain of analyticity of the scattered field in 
	$\delta$ can be extended to a neighborhood of the \textit{entire} 
	real axis (up to topological obstruction), but also that this field
	is \textit{jointly} analytic in parametric and spatial variables 
	provided that $\tepsilon(x,y,z)$ is spatially analytic. In our
	subsequent paper \cite{NichollsVo23} we extended a subset of
	these results to the three dimensional vector time--harmonic
	Maxwell equations, in particular, that the scattered field is analytic
	as a function of $\delta$ and \textit{jointly} analytic in both 
	parametric and spatial variables if $\tepsilon(x,y,z)$ is spatially 
	analytic. In the current contribution we take up not only the issue 
	of analytic continuation to perturbations $\delta$ of arbitrary (real)
	size, but also their joint analyticity with respect to spatial variables,
	which completes the analysis of these HOPE methods as applied to the
	Maxwell equations. As we shall see, this requires a significant
	enhancement of the existing technology to address variable
 coefficient 
	Maxwell equations along the lines of that presented in the work
	of Bao \& Li \cite{BaoLi22}. We note that the relevant elliptic
	estimate is a highly nontrivial generalization of that appearing
	in this latter work as it accounts for inhomogeneous terms which they
	did not analyze.
	
	The rest of the paper is organized as follows. In \S~\ref{Sec:Govern}
	we recall the governing equations and discuss transparent boundary
	conditions in \S~\ref{Sec:TransBC}. We describe the HOPE algorithm
	in \S~\ref{Sec:HOPE} and begin our theoretical developments with a statement
	of the relevant function spaces	in \S~\ref{Sec:Func}. We state and prove
	our results on parametric analytic continuation in \S~\ref{Sec:AnalCont}
	and extend these to joint analyticity in \S~\ref{Sec:JointAnal}. The crucial
	elliptic estimates upon which these results rely are established in 
	\ref{Sec:EllEst:Proof} and \ref{Sec:JointAnal:Proof}.
	
	%
	%
	
	\section{Governing Equations}
	\label{Sec:Govern}
	
	We consider materials whose electromagnetic response is modeled
	by the time--harmonic Maxwell 
	equations in three dimensions with a constant permeability
	$\mu = \mu_0$ and no currents or sources,
	\begin{gather}
		\Curl{E} - i \omega \mu_0 H = 0,
		\quad
		\Curl{H} + i \omega \eps E = 0,
		\notag \\
		\Div{\eps E} = 0,
		\quad
		\Div{H} = 0,
		\label{Eqn:Maxwell:TimeHarmonic}
	\end{gather}
	where $(E, H)$ are the electric and magnetic vector fields, and we have
	factored out time dependence of the form $\exp(-i \omega t)$ \cite{BaoLi22}.
	The permittivity $\epsilon(x,y,z)$ is biperiodic with periods $d_x$ and
	$d_y$, and is specified by
	\bes
	\eps(x,y,x) = \begin{cases}
		\epsu, & z > h, \\
		\epsv(x,y,z), & -h < z < h, \\
		\epsw, & z < -h,
	\end{cases}
	\ees
	where $\epsu, \epsw \in \Real^+$, and $\epsv(x+d_x,y+d_y,z) = \epsv(x,y,z)$,
	and
	\bes
	\lim_{z \rightarrow h-} \epsv(x,y,z) = \epsu,
	\quad
	\lim_{z \rightarrow (-h)+} \epsv(x,y,z) = \epsw.
	\ees
	Using the permittivity of vacuum, $\epsilon_0$, we can define
	\bes
	k^2_0 = \omega^2 \epsilon_0 \mu_0 = \frac{\omega^2}{c_0^2},
	\quad
	(k^m)^2 = \epsm k^2_0,
	\quad
	m \in \{u, w\},
	\ees
	and $c_0 = 1/\sqrt{\epsilon_0 \mu_0}$ is the speed of light in vacuum.
	
	This structure is illuminated from above by plane--wave incident 
	radiation of the form 
	\begin{align*}
		E^{\text{inc}}(x,y,z) & = A \exp(i \alpha x + i \beta y - i \gammau z), \\
		H^{\text{inc}}(x,y,z) & = B \exp(i \alpha x + i \beta y - i \gammau z), 
	\end{align*}
	where
	\bes
	A \cdot \kappa = 0,
	\quad
	B = \frac{1}{\omega \mu_0} \kappa \times A,
	\quad
	\Abs{A} = \Abs{B} = 1,
	\ees
	and
	\bes
	\kappa = \begin{pmatrix} \alpha \\ \beta \\ -\gammau \end{pmatrix}
	= \ku \begin{pmatrix} \sin(\theta) \cos(\phi) \\
		\sin(\theta) \sin(\phi) \\ -\cos(\theta) \end{pmatrix},
	\ees
	where $(\theta,\phi)$ are the angles of incidence.
	
	%
	%
	
	\section{Transparent Boundary Conditions}
	\label{Sec:TransBC}
	
	Following the lead of Bao \& Li \cite{BaoLi22} we use Transparent
	Boundary Conditions at $z = \pm h$ to both rigorously
	specify the appropriate far--field boundary conditions, and reduce
	the infinite domain to one of finite size. These are specified with
	Dirichlet--Neumann Operators (DNOs) which map the
	\textit{tangential} traces of the scattered electric fields at
	$z = \pm h$ to the traces of the scattered magnetic fields at 
	$z = \pm h$. Such operators are commonly called 
	Capacity Operators \cite{BaoLi22}.
	
	To summarize the developments of \cite{BaoLi22} (\S~3.2.2) we use
	the fact that $H = \frac{1}{i \omega \mu_0} \Curl{E}$ and define
	\bes
	T_u: U \rightarrow \tU,
	\quad
	T_w: W \rightarrow \tW,
	\ees
	where, for $N_u = (0,0,1)^T$ and $N_w = (0,0,-1)^T$,
	\begin{gather*}
		U := N_u \times ( \left. \Es \right|_{z=h} \times N_u),
		\quad
		\tU := \frac{1}{i \omega \mu_0} 
		( \left. \Curl{\Es} \right|_{z=h} \times N_u),
		\\
		W := N_w \times ( \left. \Es \right|_{z=-h} \times N_w),
		\quad
		\tW := \frac{1}{i \omega \mu_0}
		( \left. \Curl{\Es} \right|_{z=-h} \times N_w).
	\end{gather*}
	Using the facts that
	\begin{gather*}
		E = \Es + \Ei, \quad H = \Hs + \Hi, \quad z > h, \\
		E = \Es, \quad H = \Hs, \quad z < -h,
	\end{gather*}
	and multiplying the definitions of $\{ T_u, T_w \}$ by
	$-(i \omega \mu_0)$, we specify the Transparent Boundary Conditions
	\begin{align*}
		& \Curl{(E-\Ei)} \times N_u 
		- (i \omega \mu_0) T_u[N_u \times ((E-\Ei) \times N_u)] = 0,
		&& z = h, \\
		& \Curl{E} \times N_w 
		- (i \omega \mu_0) T_w[N_w \times (E \times N_w)] = 0,
		&& z = -h,
	\end{align*}
	or
	\begin{align*}
		& \Curl{E} \times N_u 
		- (i \omega \mu_0) T_u[N_u \times (E \times N_u)] = \phi,
		&& z = h, \\
		& \Curl{E} \times N_w 
		- (i \omega \mu_0) T_w[N_w \times (E \times N_w)] = 0,
		&& z = -h,
	\end{align*}
	where
	\be
	\label{Eqn:phi:Def}
	\phi = \Curl{\Ei} \times N_u 
	- i \omega \mu_0 \{ N_u \times (T_u[\Ei] \times N_u) \},
	\quad
	z = h.
	\ee
	
	To find a formula for $T_u$ we note that, in the upper domain 
	$\{ z > h \}$, separation of variables demands that upward
	propagating $(\alpha,\beta)$--quasiperiodic solutions of the
	Maxwell equations are
	\bes
	\Es = \sump \sumq \hat{u}_{p,q} 
	\exp(i \alpha_p x + i \beta_q y + i \gammau_{p,q} (z-h)),
	\quad
	\hat{u}_{p,q} = \begin{pmatrix} \hat{u}^x_{p,q} \\
		\hat{u}^y_{p,q} \\ \hat{u}^z_{p,q} \end{pmatrix},
	\ees
	\cite{Petit80,Yeh05} where
	\begin{gather*}
		\alpha_p = \alpha + (2 \pi/d_x) p,
		\quad
		\beta_q = \beta + (2 \pi/d_y) q,
		\\
		( \gammam_{p,q} )^2 = \epsm k_0^2 - \alpha_p^2 - \beta_q^2,
		\quad
		\ImagPart{\gammam_{p,q}} \geq 0,
		\quad
		m \in \{ u, w \}.
	\end{gather*}
	In particular, for a dielectric ($\epsm \in \Real^+$) we have
	\bes
	\gammam_{p,q} := \begin{cases}
		\sqrt{\epsm k_0^2 - \alpha_p^2 - \beta_q^2},
		& \alpha_p^2 + \beta_q^2 \leq \epsm k_0^2, \\
		i \sqrt{\alpha_p^2 + \beta_q^2 - \epsm k_0^2},
		& \alpha_p^2 + \beta_q^2 > \epsm k_0^2.
	\end{cases}
	\ees
	Using the relation
	\bes
	\begin{pmatrix} U^x(x,y) \\ U^y(x,y) \\ 0 \end{pmatrix}
	= U(x,y)
	= N_u \times (\Es(x,y,h) \times N_u)
	= \begin{pmatrix} \Esx(x,y,h) \\ \Esy(x,y,h) \\ 0 \end{pmatrix},
	\ees
	we find that
	\bes
	\hat{u}^x_{p,q} = \hat{U}^x_{p,q},
	\quad
	\hat{u}^y_{p,q} = \hat{U}^y_{p,q}.
	\ees
	To resolve $\hat{u}^z_{p,q}$ we use the divergence--free condition
	in the upper layer to deduce that
	\bes
	(i \alpha_p) \hat{u}^x_{p,q} + (i \beta_q) \hat{u}^y_{p,q}
	+ (i \gammau_{p,q}) \hat{u}^z_{p,q} = 0.
	\ees
	We now make the assumption that we are away from Rayleigh Singularities
	(commonly referred to as Wood's Anomalies)
	\bes
	\gammau_{p,q} \neq 0,
	\quad
	\forall\ p, q \in \Integer,
	\ees
	which gives
	\bes
	\hat{u}^z_{p,q} = \frac{ -\alpha_p \hat{u}^x_{p,q} 
		- \beta_q \hat{u}^y_{p,q}}{\gammau_{p,q}}
	= \frac{ -\alpha_p \hat{U}^x_{p,q} 
		- \beta_q \hat{U}^y_{p,q}}{\gammau_{p,q}}.
	\ees
	Therefore we can express
	\bes
	\Es = \sump \sumq 
	\begin{pmatrix} \hat{U}^x_{p,q} \\
		\hat{U}^y_{p,q} \\
		\frac{ -\alpha_p \hat{U}^x_{p,q} 
			- \beta_q \hat{U}^y_{p,q}}{\gammau_{p,q}}
	\end{pmatrix}
	\exp(i \alpha_p x + i \beta_q y + i \gammau_{p,q} (z-h)).
	\ees
	Now, it is a simple matter to compute
	\begin{align*}
		\tU & = \Hs \times N_u 
		= \frac{1}{i \omega \mu_0} ( \left. \Curl{\Es}
		\right|_{z=h} \times N_u ) \\
		& = \frac{1}{i \omega \mu_0} \begin{pmatrix} 
			\pz \Esx(x,y,h) - \px \Esz(x,y,h) \\
			\pz \Esy(x,y,h) - \py \Esz(x,y,h) \\ 0 \end{pmatrix}.
	\end{align*}
	For this we observe that
	\bes
	\px \Es(x,y,h) = \sump \sumq (i \alpha_p)
	\begin{pmatrix} \hat{U}^x_{p,q} \\
		\hat{U}^y_{p,q} \\
		\frac{ -\alpha_p \hat{U}^x_{p,q} 
			- \beta_q \hat{U}^y_{p,q}}{\gammau_{p,q}}
	\end{pmatrix}
	\exp(i \alpha_p x + i \beta_q y),
	\ees
	and
	\bes
	\py \Es(x,y,h) = \sump \sumq (i \beta_q)
	\begin{pmatrix} \hat{U}^x_{p,q} \\
		\hat{U}^y_{p,q} \\
		\frac{ -\alpha_p \hat{U}^x_{p,q} 
			- \beta_q \hat{U}^y_{p,q}}{\gammau_{p,q}}
	\end{pmatrix}
	\exp(i \alpha_p x + i \beta_q y),
	\ees
	and
	\bes
	\pz \Es(x,y,h) = \sump \sumq (i \gammau_{p,q})
	\begin{pmatrix} \hat{U}^x_{p,q} \\
		\hat{U}^y_{p,q} \\
		\frac{ -\alpha_p \hat{U}^x_{p,q} 
			- \beta_q \hat{U}^y_{p,q}}{\gammau_{p,q}}
	\end{pmatrix}
	\exp(i \alpha_p x + i \beta_q y).
	\ees
	Therefore,
	\begin{multline*}
		\tU = T_u[U] \\
		= \frac{1}{i \omega \mu_0} \sump \sumq
		\begin{pmatrix}
			(i \gammau_{p,q}) \hat{U}^x_{p,q} 
			+ \frac{(i \alpha_p)}{\gammau_{p,q}}
			\left\{ \alpha_p \hat{U}^x_{p,q} 
			+ \beta_q \hat{U}^y_{p,q} \right\} \\
			(i \gammau_{p,q}) \hat{U}^y_{p,q} 
			+ \frac{(i \beta_p)}{\gammau_{p,q}}
			\left\{ \alpha_p \hat{U}^x_{p,q} 
			+ \beta_q \hat{U}^y_{p,q} \right\} \\
			0
		\end{pmatrix}
		\exp(i \alpha_p x + i \beta_q y).
	\end{multline*}
	In a similar fashion
	\begin{multline*}
		\tW = T_w[W] \\
		= \frac{1}{i \omega \mu_0} \sump \sumq
		\begin{pmatrix}
			(i \gammaw_{p,q}) \hat{W}^x_{p,q} 
			+ \frac{(i \alpha_p)}{\gammaw_{p,q}}
			\left\{ \alpha_p \hat{W}^x_{p,q} 
			+ \beta_q \hat{W}^y_{p,q} \right\} \\
			(i \gammaw_{p,q}) \hat{W}^y_{p,q} 
			+ \frac{(i \beta_p)}{\gammaw_{p,q}}
			\left\{ \alpha_p \hat{W}^x_{p,q} 
			+ \beta_q \hat{W}^y_{p,q} \right\} \\
			0
		\end{pmatrix}
		\exp(i \alpha_p x + i \beta_q y).
	\end{multline*}
	
	At this point we can state our governing equations with full rigor.
	Eliminating the magnetic field from \eqref{Eqn:Maxwell:TimeHarmonic}
	and gathering our full set of governing equations we find
	the following problem to solve.
	\bse
	\label{Eqn:Max}
	\begin{align}
		& \Curl{ \Curl{ E } } - \epsv k_0^2 E = 0,
		&& \text{in $\Omega$}, \label{Eqn:Max:a} \\
		& -\Div{ \epsv k_0^2 E } = 0,
		&& \text{in $\Omega$}, \label{Eqn:Max:b} \\
		& \Curl{E} \times N_u 
		- (i \omega \mu_0) T_u[N_u \times (E \times N_u)] = \phi,
		&& \text{at $\Gamma_u$}, \label{Eqn:Max:c} \\
		& \Curl{E} \times N_w 
		- (i \omega \mu_0) T_w[N_w \times (E \times N_w)] = 0,
		&& \text{at $\Gamma_w$}, \label{Eqn:Max:d} \\
		& E(x+d_x,y+d_y,z) = \exp(i \alpha d_x + i \beta d_y) E(x,y,z), 
		\label{Eqn:Max:e}
	\end{align}
	\ese
	where
	\begin{gather*}
		\Omega := (0,d_x) \times (0,d_y) \times (-h,h), \\
		\Gamma_u := (0,d_x) \times (0,d_y) \times \{ z = h \},
		\quad
		\Gamma_w := (0,d_x) \times (0,d_y) \times \{ z = -h \}.
	\end{gather*}
	
	\begin{remark}
		Equation \eqref{Eqn:Max:b} is, of course, a simple consequence
		of the divergence operator applied to \eqref{Eqn:Max:a}. However,
		we include it \textit{explicitly} in order to highlight its
		importance in our subsequent elliptic estimates and analyticity
		theory.
	\end{remark}
	
	%
	%
	
	\section{A High--Order Perturbation of Envelopes Method}
	\label{Sec:HOPE}
	
	In our previous work \cite{Nicholls19b,NichollsVo23} we pursued
	the solution of
	\eqref{Eqn:Max} not by a classical volumetric approach, but rather 
	by a perturbative one where we thought of our configuration as a 
	\textit{small} deviation from a simpler, constant, structure,
	\bes
	\epsv(x,y,z) = \bepsilon (1 - \delta \tepsilon(x,y,z))
	= \bepsilon - \delta (\bepsilon \tepsilon(x,y,z)),
	\ees
	where $\delta \ll 1$. We showed that, provided that
	$\tepsilon$ is sufficiently smooth, the solution depends
	analytically on $\delta$ and can be expressed as a convergent
	Taylor series. Now, it is known that the coefficients of this
	series determine  the solution throughout the \textit{entire}
	domain of analyticity of $E$ which, we now show, is much larger
	than the disk of convergence of this series. In fact, it contains
	the \textit{whole} real line.
	
	To investigate this claim we study perturbations of the form
	\bes
	\epsv(x,y,z) = \bepsilon (1 - \rho \tepsilon(x,y,z)),
	\quad
	\rho \in \Real,
	\ees
	by setting $\rho = \rho_0 + \delta$ where $\rho_0$ is arbitrary,
	but real, while $\delta \ll 1$. We then write
	\begin{align*}
		\epsv(x,y,z) 
		& = \bepsilon (1 - \rho \tepsilon(x,y,z)) \\
		& = \bepsilon (1 - \rho_0 \tepsilon(x,y,z) )
		- \delta (\bepsilon \tepsilon(x,y,z)) \\
		& = \bepsilon_0(x,y,z) - \delta (\bepsilon \tepsilon(x,y,z)),
	\end{align*}
	which defines the base permittivity envelope
	\bes
	\bepsilon_0(x,y,z) := \bepsilon( 1 - \rho_0 \tepsilon(x,y,z)).
	\ees
	In contrast to our previous work \cite{NichollsVo23}, this base
	value is \textit{not} constant, but rather dependent upon all
	of the spatial variables, $(x,y,z)$.
	
	We posit that the field $E = E(x,y,z;\delta)$ depends 
	analytically upon $\delta$ so that
	\be
	\label{Eqn:EExp}
	E = E(x,y,z;\delta) 
	= \sum_{\ell=0}^{\infty} E_{\ell}(x,y,z) \delta^{\ell},
	\ee
	converges strongly in a function space. In this way we establish
	that the field, $E$, is analytic in a disk about an 
	\textit{arbitrary} real value of $\rho = \rho_0 \in \Real$.
	It is not difficult to see that these $E_{\ell}$ must satisfy
	\bse
	\label{Eqn:Max:ell}
	\begin{align}
		& \Curl{ \Curl{ E_{\ell} } } - \bepsilon_0 k_0^2 E_{\ell} 
		= \bepsilon_0 k_0^2 F_{\ell},
		&& \text{in $\Omega$}, \label{Eqn:Max:ell:a} \\
		& -\Div{ \bepsilon_0 k_0^2 E_{\ell} }
		= \Div{ \bepsilon_0 k_0^2 F_{\ell} },
		&& \text{in $\Omega$}, \label{Eqn:Max:ell:b} \\
		& \Curl{E_{\ell}} \times N_u 
		- (i \omega \mu_0) T_u[N_u \times (E_{\ell} \times N_u)] 
		= \delta_{\ell,0} \phi,
		&& \text{at $\Gamma_u$}, \label{Eqn:Max:ell:c} \\
		& \Curl{E_{\ell}} \times N_w 
		- (i \omega \mu_0) T_w[N_w \times (E_{\ell} \times N_w)] = 0,
		&& \text{at $\Gamma_w$}, \label{Eqn:Max:ell:d} \\
		& E_{\ell}(x+d_x,y+d_y,z) 
		= \exp(i \alpha d_x + i \beta d_y) E_{\ell}(x,y,z), 
		\label{Eqn:Max:ell:e}
	\end{align}
	where
	\be
	F_{\ell} 
	= F_{\ell}(x,y,z) 
	= -\frac{\bepsilon \tepsilon(x,y,z)}{\bepsilon_0(x,y,z)}
	E_{\ell-1}(x,y,z),
	\ee
	\ese
	and $\delta_{\ell,0}$ is the Kronecker delta function.
	
\begin{figure}[hbt]
		\begin{center}
			\includegraphics[width=0.45\textwidth]{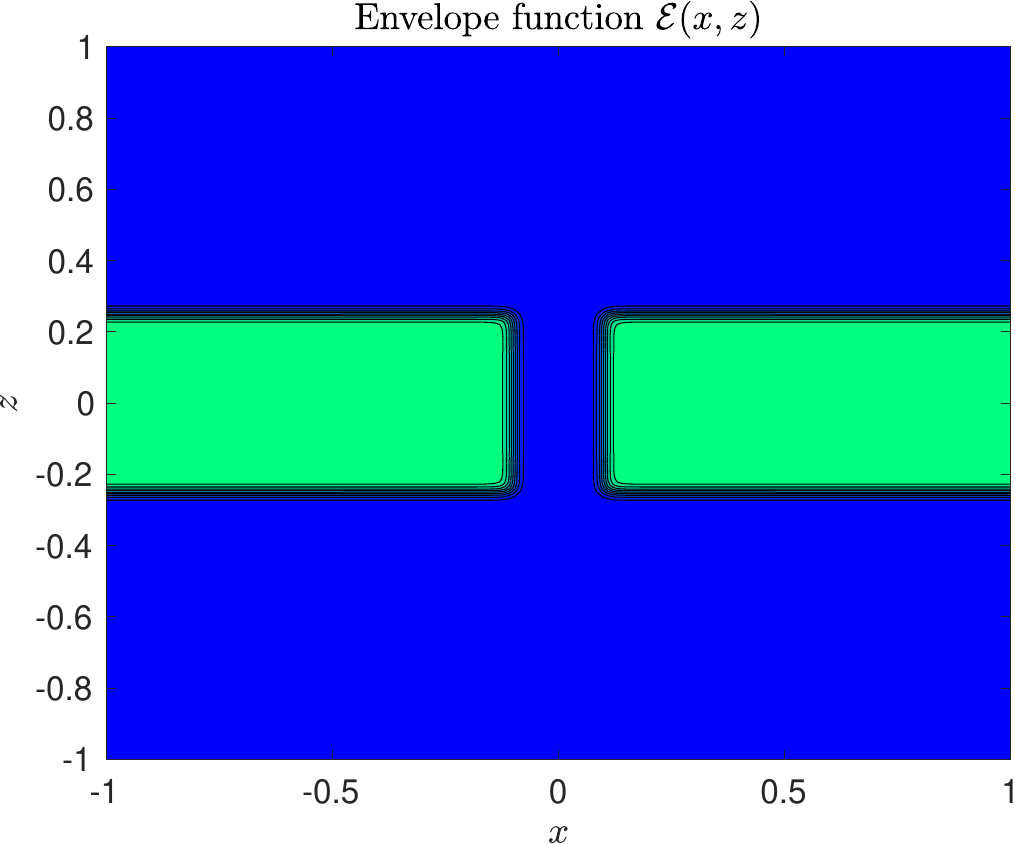}
			\quad
			\includegraphics[width=0.45\textwidth]{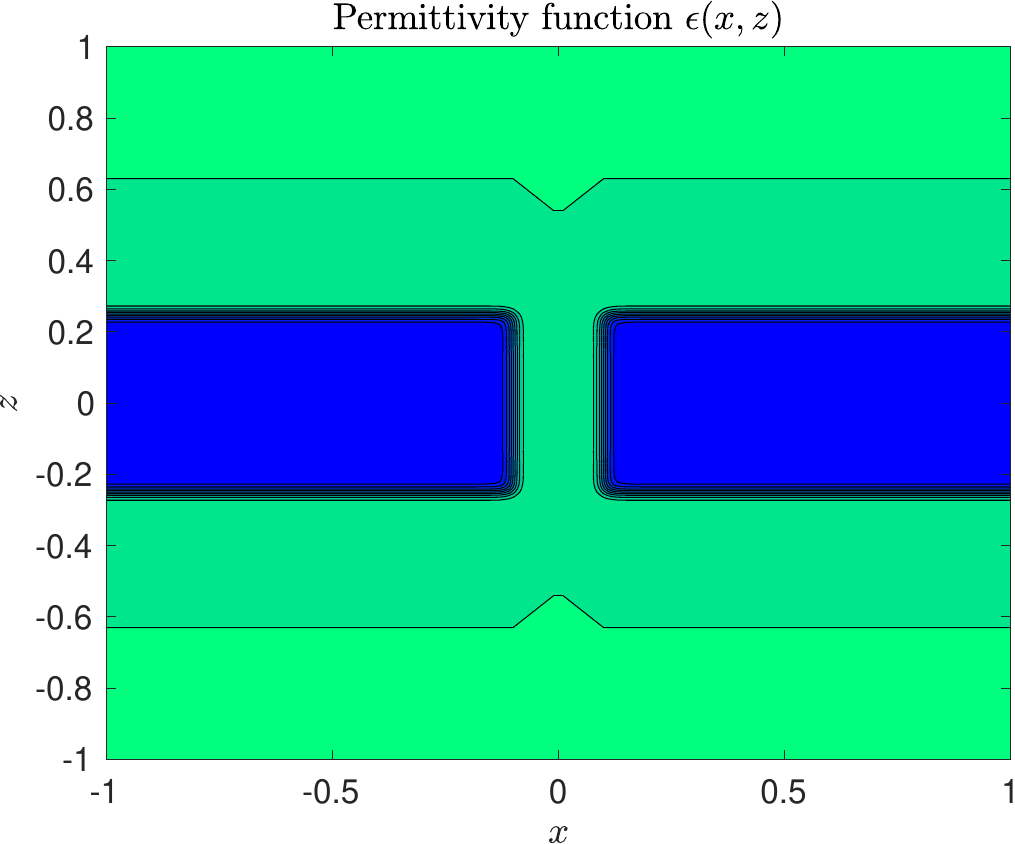}
			\caption{Contour plots of $\tepsilon(x,z)$ (left) and 
				$\epsv(x,z)$ (right).}
			\label{Fig:E:Perm}
		\end{center}
	\end{figure}

	There are many possibilities for the envelope function
	$\tepsilon(x,y,z)$ and each leads to a slightly different
	perturbation approach. For instance, consider the function
	\bes
	\Phi_{a,b}(z) := \frac{\tanh(w (z-a)) - \tanh(w (z-b))}{2},
	\ees
	with sharpness parameter $w$, which is effectively zero outside
	the interval $(a,b)$ while being essentially one inside
	$(a,b)$, c.f. \cite{Nicholls19b}. We can 
	approximate a slab of material (of permittivity $\epsilon'$)
	with thickness $2d$ and a gap of
	width $2g$ in vacuum by selecting \cite{NichollsVo23}
	\bes
	\bepsilon = 1,
	\quad
	\tepsilon(x,y,z) = \left( \frac{\bepsilon-\epsilon'}{\bepsilon} \right)
	\Phi_{-d,d}(z) \left\{ 1 - \Phi_{-g,g}(x) \right\},
	\quad
	\rho_0 = 0,
	\quad
	\delta = 1.
	\ees
	See Figure~\ref{Fig:E:Perm} with the choices $d=1/4$, $g=1/10$,
	and $w=50$ on the cell $[-1,1] \times [-1,1]$.
	%
	%
	
	%
	%
	
	\section{Function spaces}
	\label{Sec:Func}
	
	In this section, we present function spaces and theoretical
	notions that are necessary for our analysis. Due to the very
	weak formulation of the Maxwell equations we employ, espoused
	by Bao \& Li \cite{BaoLi22}, our function spaces are
	quite different from those used in our previous work
	\cite{Nicholls19b,NichollsVo23}. In fact we move to the 
	functional framework
	outlined in Section~3.3.1 of \cite{BaoLi22} of 
	$(\alpha,\beta$)--quasiperiodic $\Hcurl$ and $\Hdiv$ functions.
	More specifically
	\begin{align*}
		\Hcurl = H(\curl,\Omega)
		& = \left\{ u \in L^2(\Omega)^3\ |\ \Curl{u} \in L^2(\Omega)^3,
		\right. \\ & \quad \left.
		e^{i \alpha d_x} u(0,y,z) \times \hat{n}_x 
		= u(d_x,y,z) \times \hat{n}_x,
		\right. \\ & \quad \left.
		e^{i \beta d_y} u(x,0,z) \times \hat{n}_y 
		= u(x,d_y,z) \times \hat{n}_y
		\right\},
	\end{align*}
	where $\hat{n}_x = (1,0,0)$ and $\hat{n}_y = (0,1,0)$, and
	\bes
	\HcurlNorm{u}^2 := \Norm{u}{L^2}^2 + \Norm{\Curl{u}}{L^2}^2.
	\ees
	Additionally,
	\begin{align*}
		\Hdiv = H(\div,\Omega)
		& = \left\{ u \in L^2(\Omega)^3\ |\ \Div{u} \in L^2(\Omega),
		\right. \\ & \quad \left.
		e^{i \alpha d_x} u(0,y,z) \times \hat{n}_x 
		= u(d_x,y,z) \times \hat{n}_x,
		\right. \\ & \quad \left.
		e^{i \beta d_y} u(x,0,z) \times \hat{n}_y 
		= u(x,d_y,z) \times \hat{n}_y
		\right\},
	\end{align*}
	and
	\bes
	\HdivNorm{u}^2 := \Norm{u}{L^2}^2 + \Norm{\Div{u}}{L^2}^2.
	\ees
	Due to the particular structure of the inhomogeneous
 Maxwell equations, we
	can simplify the statement and proof of our theorems by introducing
	the $(\bepsilon_0(x,y,z) k_0^2)$--dependent space
	\bes
	X = X(\bepsilon_0,\Omega) := \left\{ u \in H(\curl,\Omega)\ |\
	(\bepsilon_0 k_0^2 u) \in H(\div,\Omega) \right\},
	\ees
	with norm
	\bes
	\XNorm{u}^2 := \HcurlNorm{u}^2 + \HdivNorm{\bepsilon_0 k_0^2 u}^2.
	\ees
	In addition, we require interfacial versions of the spaces
	$\Hcurl$ and $\Hdiv$ at the
	artificial boundaries $\Gamma_u = \{z = h\}$ and 
	$\Gamma_w = \{ z = -h \}$, namely
	\begin{align*}
		\Hmhcurl & = H^{-1/2}(\curl,\Gamma_m) \\
		& = \left\{ u \in H^{-1/2}(\Gamma_m)^3, 
		\curl_{\Gamma_m} u \in H^{-1/2}(\Gamma_m)^3,
		u^z = 0 \right\}, \\
		\Hmhdiv & = H^{-1/2}(\div,\Gamma_m) \\
		& = \left\{ u \in H^{-1/2}(\Gamma_m)^3, 
		\div_{\Gamma_m} u \in H^{-1/2}(\Gamma_m),
		u^z = 0 \right\}.
	\end{align*}
	For these, the norms can be computed \cite{BaoLi22} from
	\begin{align*}
		\HmhcurlNorm{u}^2 & := d_x d_y \sump \sumq 
		\frac{\Abs{\hat{u}^x_{p,q}}^2 + \Abs{\hat{u}^y_{p,q}}^2
			+ \Abs{ \alpha_p \hat{u}^y_{p,q} - \beta_q \hat{u}^x_{p,q} }^2}
		{\sqrt{1 + \alpha_p^2 + \beta_q^2}}, \\
		\HmhdivNorm{u}^2 & := d_x d_y \sump \sumq 
		\frac{\Abs{\hat{u}^x_{p,q}}^2 + \Abs{\hat{u}^y_{p,q}}^2
			+ \Abs{ \alpha_p \hat{u}^x_{p,q} + \beta_q \hat{u}^y_{p,q} }^2}
		{\sqrt{1 + \alpha_p^2 + \beta_q^2}}.
	\end{align*}
	We also recall the space of $s$-times continuously differentiable
	functions with H\"older norm
	\bes
	\HolderNorm{v}{s} := \max_{0 \leq \ell + r \leq s} \max_{m \in \{x, y, z\}}
	\SupNorm{\partial_x^{\ell} \partial_y^r v^m}.
	\ees
	We close with an essential result \cite{Evans10,NichollsReitich99} required
	for our later proofs.
	\begin{lemma}
		\label{Lemma_Cm}
		Let $g \in C^1(\Omega)$, $u \in H(\curl,\Omega)$, $v \in H(\div,\Omega)$,
		where $\Omega$ is a subset of $\Real^3$, then $g u \in H(\curl,\Omega)$
		and $g v \in H(\div,\Omega)$. Furthermore,
		\begin{align*}
			\HcurlNorm{g u} & \leq M(\Omega) \HolderNorm{g}{1} \HcurlNorm{u}, \\
			\HdivNorm{g v} & \leq M(\Omega) \HolderNorm{g}{1} \HdivNorm{v},
		\end{align*}
		where $M$ is some positive constant.
	\end{lemma}
	
	Finally, we recall the following elementary result
	\cite{NichollsReitich00b,Nicholls19b}.
	\begin{lemma}
		\label{Lemma:S}
		Let $s\geq0$ be an integer, then there exists a constant $S>0$ such that
		\bes
		\sum_{j=0}^s \frac{(s+1)^2}{(s-j+1)^2(j+1)^2} < S,
		\quad
		\sum_{j=0}^s \sum_{r=0}^j \frac{(s+1)^2}{(s-j+1)^2 (j-r+1)^2 (r+1)^2} < S^2.
		\ees
	\end{lemma}
	
	%
	%
	
	\section{Analytic Continuation}
	\label{Sec:AnalCont}
	
	At this point we are in a position to establish analyticity
	of the full electric field in a neighborhood of \textit{any}
	real value $\rho_0$ by demonstrating the
	analytic dependence of $E = E(x,y,z;\delta)$ upon $\delta$
	sufficiently small. More specifically, we show that the expansion 
	\eqref{Eqn:EExp} converges strongly in an appropriate function space.
	
	For this we require an elliptic estimate for our inductive proof 
	which is established in \ref{Sec:EllEst:Proof}. For future 
	convenience, we define the following differential operators associated
	to the Maxwell system
	\begin{align*}
		\cL_0 E & := \Curl{ \Curl{ E } } - \bepsilon_0(x,y,z) k_0^2 E, 
		&& \text{in $\Omega$}, \\
		\cB_m E & := \Curl{E} \times N_m 
		- (i \omega \mu_0) T_m[ N_m \times (E \times N_m) ],
		&& \text{at $\Gamma_m$},
	\end{align*}
	for $m \in \{ u, w \}$.
	As is well known \cite{BaoLi22}, the issue of \textit{uniqueness} of 
	solutions to the Maxwell problem
	\bse
	\label{Eqn:Max:Uniqueness}
	\begin{align}
		& \cL_0 V = 0, && \text{in $\Omega$}, \\
		& -\Div{\bepsilon_0 k_0^2 V} = 0, && \text{in $\Omega$}, \\
		& \cB_u V = 0, && \text{at $\Gamma_u$}, \\
		& \cB_w V  = 0, && \text{at $\Gamma_w$}, \\
		& V(x+d_x,y+d_y,z) = \exp(i \alpha d_x + i \beta d_y) V(x,y,z),
	\end{align}
	\ese
	c.f.\ \eqref{Eqn:Max:ell}, which should have only the \textit{trivial} 
	solution $V \equiv 0$, is a subtle one and certain illuminating
	frequencies $\omega$ (alternatively wavenumbers $k_0$)
	will induce non--uniqueness in some configurations.
	Unfortunately a precise characterization of the set of forbidden
	frequencies is elusive and all that is known is that it is countable
	and accumulates at infinity \cite{BaoLi22}. To accommodate this
	state of affairs we define the set of permissible configurations
	\be
	\label{Eqn:cP}
	\cP := \left\{ (\omega,\bepsilon_0)\ |\ 
	V \equiv 0\ \text{is the unique solution of
		\eqref{Eqn:Max:Uniqueness}} \right\}.
	\ee	
	With this we can now state the following fundamental elliptic
	regularity result.
	\begin{theorem}
		\label{Thm:EllEst}
		If $(\omega,\bepsilon_0) \in \cP$,
		$\bepsilon_0 \in L^{\infty}(\Omega)$,
		$(\bepsilon_0 k_0^2 F) \in H(\div,\Omega)$,
		$Q \in H^{-1/2}(\div,\Gamma_u)$, and
		$R \in H^{-1/2}(\div,\Gamma_w)$, 
		then there exists a unique solution 
  $E \in X(\bepsilon_0,\Omega)$ of
		\bse
		\label{Eqn:EllEst}
		\begin{align}
			& \cL_0 E = \bepsilon_0 k_0^2 F, && \text{in $\Omega$}, 
			\label{Eqn:EllEst:a} \\
			& -\Div{ \bepsilon_0 k_0^2 E } = \Div{\bepsilon_0 k_0^2 F}, 
			&& \text{in $\Omega$}, \label{Eqn:EllEst:b} \\
			& \cB_u E = Q, && \text{at $\Gamma_u$}, \\
			& \cB_w E = R, && \text{at $\Gamma_w$}, \\
			& E(x+d_x,y+d_y,z) = \exp(i \alpha d_x + i \beta d_y) E(x,y,z),
		\end{align}
		\ese
		satisfying
		\be
		\label{Eqn:EllEst:Est}
		\XNorm{E} \leq C_e \left( \HdivNorm{\bepsilon_0 k_0^2 F}
		+ \HmhdivNorm{Q} + \HmhdivNorm{R} \right),
		\ee
		where $C_e > 0$ is a positive constant.
	\end{theorem}
	
	We can now prove the following analytic continuation result.
	\begin{theorem}
		\label{Thm:AnalCont}
		If $(\omega,\bepsilon_0) \in \cP$, 
		$\bepsilon_0 \in L^{\infty}(\Omega)$,
		and $(\tepsilon/\bepsilon_0) \in C^1(\Omega)$
		then the series \eqref{Eqn:EExp} converges strongly. More precisely,
		\be
		\label{Eqn:AnalCont}
		\XNorm{E_{\ell}} \leq K B^{\ell},
		\quad \forall\ \ell \geq 0,
		\ee
		for some constants $K, B >0$.
	\end{theorem}
	\begin{proof}
		We prove the estimate \eqref{Eqn:AnalCont} by induction. 
		For $\ell = 0$ the system \eqref{Eqn:Max:ell} can be written as
		\begin{align*}
			& \cL_0 E_0 = 0, && \text{in $\Omega$}, \\
			& -\Div{ \bepsilon_0 k_0^2 E_0 } = 0, && \text{in $\Omega$}, \\
			& \cB_u E_0 = \phi, && \text{at $\Gamma_u$}, \\
			& \cB_w E_0 = 0, && \text{at $\Gamma_w$}, \\
			& E_0(x+d_x,y+d_y,z) = \exp(i \alpha d_x + i \beta d_y) E_0(x,y,z).
		\end{align*}
		We can apply Theorem~\ref{Thm:EllEst} with $F \equiv 0$,
		$Q = \phi$, and $R \equiv 0$ to obtain \eqref{Eqn:AnalCont}
		\bes
		\XNorm{E_{\ell}} \leq C_e \HmhdivNorm{\phi} =: K.
		\ees
		
		Next, we assume that \eqref{Eqn:AnalCont} is true for all $\ell < L$
		and apply Theorem~\ref{Thm:EllEst} to the system \eqref{Eqn:Max:ell}
		for $E_L$ with $F = -\bepsilon (\tepsilon/\bepsilon_0) E_{L-1}$ and
		$Q \equiv R \equiv 0$. This gives
		\begin{align*}
			\XNorm{E_L}
			& \leq C_e \HdivNorm{\bepsilon_0 k_0^2 F_L} \\
			& \leq C_e \HdivNorm{(\bepsilon_0 k_0^2) 
				(-\bepsilon (\tepsilon/\bepsilon_0) E_{L-1})} \\
			& \leq C_e \bepsilon M \HolderNorm{\tepsilon/\bepsilon_0}{1}
			\HdivNorm{\bepsilon_0 k_0^2 E_{L-1}} \\
			& \leq C_e \bepsilon M \HolderNorm{\tepsilon/\bepsilon_0}{1}
			K B^{L-1} \\
			& \leq K B^L,
		\end{align*}
		provided that
		\bes
		B > C_e \bepsilon M \HolderNorm{\tepsilon/\bepsilon_0}{1},
		\ees
		and we are done.
	\end{proof}
	
	\begin{remark}
		At this point we comment on the smoothness requirements we make
		on $\tepsilon(x,y,z)$ and 
		$\bepsilon_0(x,y,z) = \bepsilon (1 - \rho_0 \tepsilon(x,y,z))$.
		First, we ask that $\bepsilon_0 \in L^{\infty}$ which imposes a
		(completely appropriate) boundedness requirement on $\tepsilon(x,y,z)$.
		However, we also ask that $\tepsilon(x,y,z)/\bepsilon_0(x,y,z) \in C^1$
		which puts more complicated size requirements on $\rho_0 \tepsilon(x,y,z)$
		measured in the $C^1$ norm. However, it is clear that
		$\Abs{\rho_0} \HolderNorm{\tepsilon(x,y,z)}{1}$ cannot be large and
		must not be too small.
	\end{remark}
	
	From this we can derive the exponential order of convergence
	of this HOPE method. More precisely, defining the $L$--th 
	partial sum of \eqref{Eqn:EExp}
	\bes
	E^L(x,y,z;\delta) := \sum_{\ell=0}^{L} E_{\ell}(x,y,z) \delta^{\ell},
	\ees
	we obtain the following error estimate.
	\begin{theorem}
		If $(\omega,\bepsilon_0) \in \cP$,
		$\bepsilon_0 \in L^{\infty}(\Omega)$,
		and $E$ is the unique solution
		of \eqref{Eqn:Max}, under the assumptions of Theorem~\ref{Thm:AnalCont}
		we have the estimate
		\bes
		\XNorm{E-E_L} \leq \tilde{K} (B \delta)^{L+1},
		\ees
		for some constants $\tilde{K}, B > 0$, provided that 
		$\Abs{\delta} < 1/B$.
	\end{theorem}
	\begin{proof}
		Since
		\bes
		E(x,y,z) - E^L(x,y,z) 
		= \sum_{\ell = L+1}^{\infty} E_{\ell}(x,y,z) \delta^{\ell},
		\ees
		we have, by Theorem~\ref{Thm:AnalCont},
		\bes
		\XNorm{E - E^L} \leq 
		\sum_{\ell = L+1}^{\infty} \XNorm{E_{\ell}} \delta^{\ell}
		\leq \sum_{\ell = L+1}^{\infty} K B^{\ell} \delta^{\ell}.
		\ees
		By gathering terms and re--indexing we have
		\bes
		\XNorm{E - E^L}
		\leq K (B \delta)^{L+1} \sum_{\ell=0}^{\infty} (B \delta)^{\ell}
		\leq \tilde{K} (B \delta)^{L+1},
		\ees
		for $\Abs{\delta B} < 1$, where $\tilde{K} := K/(1- (B \delta))$, and 
		we have used the elementary fact that
		\bes
		\sum_{\ell=0}^{\infty} \alpha^{\ell} = \frac{1}{1 - \alpha},
		\ees
		provided that $\Abs{\alpha} < 1$.
	\end{proof}
	
	%
	%
	
	\section{Joint Analyticity}
	\label{Sec:JointAnal}
	
	In this section we show that the solution $E(x,y,z;\delta)$
	of \eqref{Eqn:Max} is jointly analytic with respect to both
	parameter, $\delta$, and spatial variables, $(x,y,z)$. For this
	we need an appropriate notion of analyticity which we give in the 
	following definition of $C^{\omega}_q$.
	\begin{definition}
		Given an integer $q \geq 0$, if the functions $f = f(x,y)$ and 
		$\tepsilon = \tepsilon(x,y,z)$ are \textit{real analytic}
		and satisfy the following estimates
		\begin{align*}
			\HolderNorm{\frac{\px^r \py^t}{(r+t)!} f}{q} 
			& \leq C_f \frac{\eta^r}{(r+1)^2} \frac{\theta^t}{(t+1)^2}, \\
			\HolderNorm{\frac{\px^r \py^t \pz^s}{(r+t+s)!} \tepsilon}{q} 
			& \leq C_{\tepsilon} \frac{\eta^r}{(r+1)^2} \frac{\theta^t}{(t+1)^2}
			\frac{\zeta^s}{(s+1)^2},
		\end{align*}
		for all $r, t, s \geq 0$
		and constants $C_f, C_{\tepsilon}, \eta, \theta, \zeta > 0$,
		then $f \in C^{\omega}_q(\Gamma_m)$, $m \in \{u,w\}$,
		and $\tepsilon \in C^{\omega}_q(\Omega)$.
	\end{definition}
	
	Here $C^{\omega}_q$ is the space of real analytic functions
	with radius of convergence (specified by $\eta$, $\theta$, and
	$\zeta$) measured in the $C^q$ norm.
	It is clear that the incident radiation function $\phi$,
	\eqref{Eqn:phi:Def}, is jointly
	analytic in $x$ and $y$ as we now explicitly state.
	\begin{lemma}
		The function $\phi(x,y)$ defined in \eqref{Eqn:phi:Def}
		is real analytic and satisfies
		\bes
		\HmhdivNorm{ \frac{\px^r \py^t}{(r+t)!} \phi }
		\leq C_{\phi} \frac{\eta^r}{(r+1)^2} \frac{\theta^t}{(t+1)^2},
		\ees
		for all $r, t \geq 0$ and some constants $C_{\phi}, \eta, \theta > 0$.
	\end{lemma}
	
	Now we present the fundamental elliptic estimate which is required
	in our estimates. (It is proven in \ref{Sec:JointAnal:Proof}.)
	\begin{theorem}
		\label{Thm:JointAnal:EllEst}
		Given any integer $q \geq 0$, if
		$(\omega,\bepsilon_0) \in \cP$,
		$\bepsilon_0 \in C^{\omega}_q(\Omega)$, such that
		\bes
		\HolderNorm{\frac{\px^r \py^t \pz^s}{(r+t+s)!} \bepsilon_0}{q}
		\leq C_{\bepsilon_0} \frac{\eta^r}{(r+1)^2} \frac{\theta^t}{(t+1)^2}
		\frac{\zeta^s}{(s+1)^2},
		\ees
		for all $r, t, s \geq 0$ and some constants 
		$C_{\bepsilon_0}, \eta, \theta, \zeta > 0$,
		and $(\bepsilon_0 k_0^2 F) \in C^{\omega}(\Omega)$ such that
		\bes
		\HdivNorm{\frac{\px^r \py^t \pz^s}{(r+t+s)!} 
			\left[ \bepsilon_0 k_0^2 F \right]}
		\leq C_F \frac{\eta^r}{(r+1)^2} \frac{\theta^t}{(t+1)^2}
		\frac{\zeta^s}{(s+1)^2},
		\ees
		for all $r, t, s \geq 0$ and some constant $C_F > 0$, and
		$Q \in C^{\omega}(\Gamma_u)$ and
		$R \in C^{\omega}(\Gamma_w)$ satisfying 
		\begin{align*}
			& \HmhdivNorm{\frac{\px^r \py^t}{(r+t)!} Q}
			\leq C_Q \frac{\eta^r}{(r+1)^2} \frac{\theta^t}{(t+1)^2}, \\
			& \HmhdivNorm{\frac{\px^r \py^t}{(r+t)!} R}
			\leq C_R \frac{\eta^r}{(r+1)^2} \frac{\theta^t}{(t+1)^2},
		\end{align*}
		for all $r, t \geq 0$ and some constants $C_R, C_Q > 0$.
		Then, there exists a unique solution
		$E \in C^{\omega}(\Omega)$ of 
		\bse
		\label{Eqn:JointAnal:EllEst}
		\begin{align}
			& \cL_0 E = \bepsilon_0 k_0^2 F,
			&& \text{in $\Omega$}, \\
			& -\Div{\bepsilon_0 k_0^2 E} = \Div{ \bepsilon_0 k_0^2 F}, 
			&& \text{in $\Omega$}, \\
			& \cB_u E = Q, && \text{at $\Gamma_u$}, \\
			& \cB_w E = R, && \text{at $\Gamma_w$}, \\
			& E(x+d_x,y+d_y, z) = e^{i \alpha d_x + i \beta d_y} E(x,y,z),
		\end{align}
		\ese
		satisfying 
		\be
		\label{Eqn:JointAnal:EllEst:Estimate}
		\XNorm{ \frac{\px^r \py^t \pz^s}{(r+t+s)!} E } 
		\leq \underline{C}_e \frac{\eta^r}{(r+1)^2}
		\frac{\theta^t}{(t+1)^2} \frac{\zeta^s}{(s+1)^2},
		\ee
		for all $r, t, s \geq 0$ where
		\bes
		\underline{C}_e = \overline{C} (C_F + C_Q + C_R) > 0,
		\ees
		and $\overline{C} > 0$ is a constant.
	\end{theorem}
	
	We now give the recursive estimate which is essential for our
	joint analyticity result.
	\begin{lemma}
		\label{Lemma:JointAnal:IndEst}
		Given any integer $q \geq 0$, if
		$(\omega,\bepsilon_0) \in \cP$;
		$\bepsilon_0, (\tepsilon/\bepsilon_0) \in C^{\omega}_q(\Omega)$,
		such that
		\begin{align*}
			\HolderNorm{\frac{\px^r \py^t \pz^s}{(r+t+s)!} \bepsilon_0}{q}
			& \leq C_{\bepsilon_0} \frac{\eta^r}{(r+1)^2} \frac{\theta^t}{(t+1)^2}
			\frac{\zeta^s}{(s+1)^2}, \\
			\HolderNorm{\frac{\px^r \py^t \pz^s}{(r+t+s)!} 
				[\tepsilon/\bepsilon_0]}{q}
			& \leq C_{\tepsilon/\bepsilon_0} 
			\frac{\eta^r}{(r+1)^2} \frac{\theta^t}{(t+1)^2}
			\frac{\zeta^s}{(s+1)^2},
		\end{align*}
		for all $r, t, s \geq 0$ and some constants 
		$C_{\bepsilon_0} C_{\tepsilon/\bepsilon_0}, \eta, \theta, \zeta > 0$,
		and 
		\bes
		\XNorm{ \frac{\px^r \py^t \pz^s}{(r+t+s)!} E_{\ell} }
		\leq K B^{\ell} \frac{\eta^r} {(r+1)^2} \frac{\theta^t}{(t+1)^2}
		\frac{\zeta^s}{(s+1)^2},
		\quad \forall\ \ell < L,
		\ees
		for all $r, t, s \geq 0$ and for some constants $K, B > 0$. Then, 
		\be
		\label{Eqn:Recur}
		\HdivNorm{ \frac{\px^r \py^t \pz^s}{(r+t+s)!} 
			\left[ (\bepsilon_0 k_0^2) F_L \right] }
		\leq \tilde{C} K B^{L-1} \frac{\eta^r}{(r+1)^2}
		\frac{\theta^t}{(t+1)^2} \frac{\zeta^s}{(s+1)^2},
		\ee
		for all $r, t, s \geq 0$ and some constant $\tilde{C}>0$.
	\end{lemma}
	\begin{proof}
		Using Leibniz's rule, we have that
		\begin{align*}
			\frac{\px^r \py^t \pz^{s}}{(r+t+s)!} 
			& \left[ (\bepsilon_0 k_0^2) F_L \right]
			= \frac{\px^r \py^t \pz^{s}}{(r+t+s)!} 
			\left[ (\bepsilon_0 k_0^2) (-\bepsilon (\tepsilon/\bepsilon_0) 
			E_{L-1}) \right] \\
			& = -\bepsilon \frac{\px^r \py^t \pz^{s}}{(r+t+s)!} 
			\left[ (\tepsilon/\bepsilon_0) k_0^2 \bepsilon_0 E_{L-1} \right] \\
			& = -\bepsilon 
			\frac{r! t! s!}{(r+t+s)!} \sum_{j=0}^{r} \sum_{k=0}^{t}
			\sum_{\ell=0}^{s} \left( \frac{\px^{r-j}}{(r-j)!}
			\frac{\py^{t-k}}{(t-k)!} \frac{\pz^{s-\ell}}{(s-\ell)!}
			\left[ \tepsilon/\bepsilon_0 \right]
			\right) \\
			& \quad \times \left( \frac{\px^j}{j!} \frac{\py^k}{k!} 
			\frac{\pz^{\ell}}{\ell!} \left[ \bepsilon_0 k_0^2 E_{L-1}
			\right] \right).
		\end{align*}
		Using the inequality $r! t! s! \leq (r+t+s)!$, we obtain
		\begin{multline*}
			\HdivNorm{\frac{\px^r \py^t \pz^{s}}{(r+t+s)!} 
				\left[ (\bepsilon_0 k_0^2) F_L \right]} \\
			\leq \bepsilon
			\sum_{j=0}^{r} \sum_{k=0}^{t} \sum_{\ell=0}^{s}
			\HdivNorm{\left( \frac{\px^{r-j}}{(r-j)!} \frac{\py^{t-k}}{(t-k)!}
				\frac{\pz^{s-\ell}}{(s-\ell)^2} \left[ \tepsilon/\bepsilon_0 \right] \right)
				\left( \frac{\px^j}{j!} \frac{\py^k}{k!} \frac{\pz^{\ell}}{\ell!}
				\left[ \bepsilon_0 k_0^2 E_{L-1} \right] \right)}
		\end{multline*}
		Continuing
		\begin{align*}
			\HdivNorm{\frac{\px^r \py^t \pz^{s}}{(r+t+s)!} 
				\left[ (\bepsilon_0 k_0^2) F_L \right]}
			& \leq \bepsilon M 
			\sum_{j=0}^{r} \sum_{k=0}^{t} \sum_{\ell=0}^{s}
			\HolderNorm{\frac{\px^{r-j}}{(r-j)!} \frac{\py^{t-k}}{(t-k)!}
				\frac{\pz^{s-\ell}}{(s-\ell)^2} \left[ \tepsilon/\bepsilon_0 \right]}{1} \\
			& \quad \times
			\HdivNorm{\frac{\px^j}{j!} \frac{\py^k}{k!} \frac{\pz^{\ell}}{\ell!}
				\left[ \bepsilon_0 k_0^2 E_{L-1} \right]} \\
			& \leq \bepsilon M
			\sum_{j=0}^{r} \sum_{k=0}^{t} \sum_{\ell=0}^{s}
			C_{\tepsilon/\bepsilon_0} \frac{\eta^{r-j}}{(r-j+1)^2}
   \\ & \quad \times
			\frac{\theta^{t-k}}{(t-k+1)^2} \frac{\zeta^{s-\ell}}{(s-\ell+1)^2}
			\\
			& \quad \times K B^{L-1} \frac{\eta^j}{(j+1)^2}
			\frac{\theta^k}{(k+1)^2} \frac{\zeta^{\ell}}{(\ell+1)^2} \\
			& \leq \bepsilon M C_{\tepsilon/\bepsilon_0} K B^{L-1} 
			\frac{\eta^r}{(r+1)^2} \frac{\theta^t}{(t+1)^2} 
			\frac{\zeta^s}{(s+1)^2} \\
			& \quad \times \sum_{j=0}^{r} \frac{(r+1)^2}{(r-j+1)^2 (j+1)^2}
			\sum_{k=0}^{t} \frac{(t+1)^2}{(t-k+1)^2(k+1)^2} \\
			& \quad \times \sum_{\ell=0}^{s} 
			\frac{(s+1)^2}{(s-\ell+1)^2(\ell+1)^2} \\
			& \leq \bepsilon M C_{\tepsilon/\bepsilon_0} M S^3 K B^{L-1}
			\frac{\eta^r}{(r+1)^2} \frac{\theta^t}{(t+1)^2}
			\frac{\zeta^s}{(s+1)^2},
		\end{align*}
		where $S$ comes from Lemma~\ref{Lemma:S}. If we choose
		\bes
		\tilde{C} \geq \bepsilon M C_{\tepsilon/\bepsilon_0} S^3,
		\ees
		the proof is complete.
	\end{proof}
	
	We conclude with our joint analyticity theorem.
	\begin{theorem}
		Given any integer $q \geq 0$, if
		$(\omega,\bepsilon_0) \in \cP$;
		$\bepsilon_0, (\tepsilon/\bepsilon_0) \in C^{\omega}_q(\Omega)$,
		such that
		\begin{align*}
			\HolderNorm{\frac{\px^r \py^t \pz^s}{(r+t+s)!} \bepsilon_0}{q}
			& \leq C_{\bepsilon_0} \frac{\eta^r}{(r+1)^2} \frac{\theta^t}{(t+1)^2}
			\frac{\zeta^s}{(s+1)^2}, \\
			\HolderNorm{\frac{\px^r \py^t \pz^s}{(r+t+s)!} 
				[\tepsilon/\bepsilon_0]}{q}
			& \leq C_{\tepsilon/\bepsilon_0} 
			\frac{\eta^r}{(r+1)^2} \frac{\theta^t}{(t+1)^2}
			\frac{\zeta^s}{(s+1)^2},
		\end{align*}
		for all $r, t, s \geq 0$ and some constants 
		$C_{\bepsilon_0} C_{\tepsilon/\bepsilon_0}, \eta, \theta, \zeta > 0$.
		Then the series \eqref{Eqn:EExp} converges strongly. Moreover the 
		$E_{\ell}(x,y,z)$ satisfy the joint analyticity estimate
		\be
		\label{Eqn:JointAnal}
		\XNorm{ \frac{\px^r \py^t \pz^s}{(r+t+s)!} E_{\ell} }
		\leq K B^{\ell} \frac{\eta^r}{(r+1)^2} \frac{\theta^t}{(t+1)^2}
		\frac{\zeta^s}{(s+1)^2},
		\ee
		for all $\ell, r, t, s \geq 0$ and constants $K, B > 0$.
	\end{theorem}
	\begin{proof}
		We prove \eqref{Eqn:JointAnal} by induction, beginning with $\ell = 0$.
		Applying Theorem~\ref{Thm:JointAnal:EllEst} with $F \equiv 0$,
		$Q = \phi$, and $R \equiv 0$ we obtain
		\bes
		\XNorm{ \frac{\px^r \py^t \pz^s}{(r+t+s)!} E_0 }
		\leq C_{\phi} \frac{\eta^r}{(r+1)^2} \frac{\theta^t}{(t+1)^2},
		\ees
		for all $r, t, s \geq 0$ which establishes \eqref{Eqn:JointAnal}
		with $K := C_{\phi}$.
		
		Next we assume that \eqref{Eqn:JointAnal} is valid for all $\ell < L$.
		With $\ell = L$ we invoke Lemma~\ref{Lemma:JointAnal:IndEst} and 
		apply Theorem~\ref{Thm:JointAnal:EllEst} with $F \equiv F_L$,
		$C_F = \tilde{C} K B^{L-1}$, $Q \equiv 0$, and $R \equiv 0$,
		to arrive at
		\bes
		\XNorm{ \frac{\px^r \py^t \pz^s}{(r+t+s)!} E_L }
		\leq \overline{C} \tilde{C} K B^{L-1} \frac{\eta^r}{(r+1)^2}
		\frac{\theta^t}{(t+1)^2} \frac{\zeta^s}{(s+1)^2},
		\ees
		for all $r, t, s \geq 0$.
		The proof is complete by choosing $B > \overline{C} \tilde{C}$.
	\end{proof}

 %
%

\section*{Declarations}

The authors have no relevant financial or non--financial
interests to disclose.
All authors contributed to the study conception and design.
Material preparation, data collection and analysis were 
performed by all authors. The first 
draft of the manuscript was written by D.P.\ Nicholls and
all authors commented on previous versions of the manuscript.
All authors read and approved the final manuscript.
	
	%
	%
	
	\appendix
	
	%
	%
	
	\section{The Proof of Theorem~\ref{Thm:EllEst}}
	\label{Sec:EllEst:Proof}
	
	Following Bao \& Li \cite{BaoLi22} we dot the Maxwell
	equation \eqref{Eqn:EllEst:a} with 
	$\bar{w} \in \Hcurl$ and integrate over $\Omega$
	\bes
	\int_{\Omega} (\Curl{\Curl{E}}) \cdot \bar{w} \dV
	- k_0^2 \int_{\Omega} \bepsilon_0 E \cdot \bar{w} \dV 
	= \int_{\Omega} (\bepsilon_0 k_0^2 F) \cdot \bar{w} \dV.
	\ees
	We now use the first vector Green theorem \cite{BaoLi22}
	\bes
	\int_{\Omega} u \cdot ( \Curl{ \sigma \Curl{v} } ) \dV
	= \int_{\Omega} \sigma \Curl{u} \cdot \Curl{v} \dV
	- \oint_{\partial \Omega} \sigma (u \times \Curl{v}) \cdot \nu \dS,
	\ees
	to obtain (with $u=\bar{w}$, $v=E$, and $\sigma=1$)
	\begin{multline*}
		\int_{\Omega} \Curl{E} \cdot \Curl{\bar{w}} \dV
		- \oint_{\partial \Omega} (\bar{w} \times \Curl{E}) \cdot N \dS \\
		- k_0^2 \int_{\Omega} \epsv E \cdot \bar{w} \dV 
		= \int_{\Omega} (\bepsilon_0 k_0^2 F) \cdot \bar{w} \dV.
	\end{multline*}
	Using the triple product identity
	$a \cdot (b \times c) = b \cdot (c \times a)$ we find
	\begin{multline*}
		\int_{\Omega} \Curl{E} \cdot \Curl{\bar{w}} \dV
		- \oint_{\partial \Omega} \bar{w} \cdot (\Curl{E} \times N) \dS \\
		- k_0^2 \int_{\Omega} \epsv E \cdot \bar{w} \dV 
		= \int_{\Omega} (\bepsilon_0 k_0^2 F) \cdot \bar{w} \dV.
	\end{multline*}
	As $E$ and $w$ are $(\alpha,\beta)$--quasiperiodic,
 the contributions from the
	boundary of $\Omega$ reduce to $\Gamma_u$ and $\Gamma_w$,
	\begin{multline*}
		\int_{\Omega} \Curl{E} \cdot \Curl{\bar{w}} \dV
		- \int_{\Gamma_u} \bar{w} \cdot (\Curl{E} \times N_u) \dS
		- \int_{\Gamma_w} \bar{w} \cdot (\Curl{E} \times N_w) \dS \\
		- k_0^2 \int_{\Omega} \epsv E \cdot \bar{w} \dV 
		= \int_{\Omega} (\bepsilon_0 k_0^2 F) \cdot \bar{w} \dV.
	\end{multline*}
	Using the boundary conditions at $z = \pm h$ we find
	\begin{multline*}
		\int_{\Omega} \Curl{E} \cdot \Curl{\bar{w}} \dV
		- k_0^2 \int_{\Omega} \epsv E \cdot \bar{w} \dV \\
		- \int_{\Gamma_u} \bar{w} \cdot \left\{
		(i \omega \mu_0) T_u[N_u \times (E \times N_u)]
		+ Q \right\} \dS \\
		- \int_{\Gamma_w} \bar{w} \cdot \left\{
		(i \omega \mu_0) T_w[N_w \times (E \times N_w)]
		+ R \right\} \dS
		= \int_{\Omega} (\bepsilon_0 k_0^2 F) \cdot \bar{w} \dV.
	\end{multline*}
	We write this as
	\be
	a(E,w) = L[w],
	\ee
	where
	\begin{align*}
		a(E,w) 
		& = \InnerOmega{\Curl{E}}{\Curl{w}}
		- k_0^2 \InnerOmega{\epsv E}{w} \\
		& \quad
		- \InnerGammau{(i \omega \mu_0) T_u[N_u \times (E \times N_u)]}{w}
		- \InnerGammaw{(i \omega \mu_0) T_w[N_w \times (E \times N_w)]}{w},
	\end{align*}
	and
	\bes
	L[w] = \InnerOmega{(\bepsilon_0 k_0^2 F)}{w}
	+ \InnerGammau{Q}{w}
	+ \InnerGammaw{R}{w}.
	\ees
	In these we use the duality pairings
	\bes
	\InnerOmega{u}{v} := \int_{\Omega} u \cdot \bar{v} \dV,
	\quad
	\InnerGammam{u}{v} := \int_{\Gamma_m} u \cdot \bar{v} \dS.
	\ees
	
	We now seek a solution, $E \in \Hcurl$, of this weak formulation 
	by writing the $E, w \in \Hcurl$ in the form
	\begin{gather*}
		E = \vec{u} + \Grad{u},
		\quad
		w = \vec{v} + \Grad{v},
		\\
		\vec{u}, \vec{v} \in \HH,
		\quad
		u, v \in \Hzone,
		\quad
		\Grad{u}, \Grad{v} \in \HHp.
	\end{gather*}
	In these we use the spaces $\HH$ and $\HHp$ defined in 
	Bao \& Li \cite{BaoLi22}
	\begin{align*}
		\HH & = \left\{ \vec{u} \in \Hcurl\ |\
		\Div{\bepsilon_0 \vec{u}} = 0\ \text{in $\Omega$},
		\right. \\ & \quad \left.
		-k_0^2 \bepsilon_0 \vec{u} \cdot N_m 
		+ \DivGm{ (i \omega \mu_0) T_m[\vec{u}_{\Gamma_m}] } = 0\
		\text{at $\Gamma_m$} \right\}, \\
		\HHp & = \left\{ \vec{u}\ |\ \vec{u} = \Grad{u},
		u \in H_0^1 \right\} \\
		H_0^1 & = \left\{ u = H^1\ |\ \int_{\Omega} u \dV = 0 
		\right\}.
	\end{align*}
	Inserting these into our weak form we find
	\bes
	a(\vec{u}+\Grad{u},\vec{v}+\Grad{v}) = L[\vec{v}+\Grad{v}].
	\ees
	Using the fact that $a(\vec{u},\Grad{v}) = 0$, \cite{BaoLi22},
	we have
	\bes
	a(\vec{u},\vec{v}) + a(\Grad{u},\vec{v}) + a(\Grad{u},\Grad{v})
	= L[\vec{v}] + L[\Grad{v}],
	\ees
	which we solve in two phases: First for the terms involving
	$\Grad{v}$, and then for the remaining terms.
	
	%
	%
	
	\subsection{Finding a Solution I: $\HHp$}
	
	We begin by determining $u$ from
	\bes
	a(\Grad{u},\Grad{v}) = L[\Grad{v}].
	\ees
	To proceed we note the following result.
	\begin{lemma}
		\label{Lemma:IBP}
		We have
		\begin{align*}
			\InnerGammau{Q}{\Grad{v}}
			& = \InnerGammau{-\DivGu{Q}}{v}, \\
			\InnerGammaw{R}{\Grad{v}}
			& = \InnerGammaw{-\DivGw{R}}{v}, \\
			\InnerOmega{(\bepsilon_0 k_0^2 F)}{\Grad{v}}
			& = \InnerOmega{-\Div{(\bepsilon_0 k_0^2 F)}}{v} 
			+ \InnerGammau{(\bepsilon_0 k_0^2 F) \cdot N}{v}
			+ \InnerGammaw{(\bepsilon_0 k_0^2 F) \cdot N}{v}.
		\end{align*}
	\end{lemma}
	\begin{proof}
		From integration--by--parts,
		\begin{align*}
			\InnerGammau{Q}{\Grad{v}}
			& = \int_{\Gamma_u} Q \cdot (\GradGu{\bar{v}}) \dS \\
			& = \int_{\Gamma_u} \DivGu{Q \cdot \bar{v}} \dS
			- \int_{\Gamma_u} \DivGu{Q} \bar{v} \dS \\
			& = \left[ Q \cdot \bar{v} \right]_{x=0}^{x=d}
			- \int_{\Gamma_u} \DivGu{Q} \bar{v} \dS \\
			& = \InnerGammau{-\DivGu{Q}}{v},
		\end{align*}
		by the periodicity of the product $Q \cdot \bar{v}$. In a
		similar fashion we have
		\begin{align*}
			\InnerOmega{F}{\Grad{v}}
			& = \int_{\Omega} (\bepsilon_0 k_0^2 F) \cdot \Grad{\bar{v}} \dV \\
			& = \int_{\Omega} \Div{ (\bepsilon_0 k_0^2 F) \bar{v} } \dV
			- \int_{\Omega} \Div{(\bepsilon_0 k_0^2 F)} \bar{v} \dV \\
			& = \oint_{\partial \Omega} ((\bepsilon_0 k_0^2 F) \cdot N) \bar{v} \dS
			- \int_{\Omega} \Div{(\bepsilon_0 k_0^2 F)} \bar{v} \dV \\
			& = -\int_{\Omega} \Div{(\bepsilon_0 k_0^2 F)} \bar{v} \dV
			+ \int_{\Gamma_u} ((\bepsilon_0 k_0^2 F) \cdot N) \bar{v} \dS
			+ \int_{\Gamma_w} ((\bepsilon_0 k_0^2 F) \cdot N) \bar{v} \dS \\
			& = \InnerOmega{-(\bepsilon_0 k_0^2 F)}{v} 
			+ \InnerGammau{(\bepsilon_0 k_0^2 F) \cdot N}{v}
			+ \InnerGammaw{(\bepsilon_0 k_0^2 F) \cdot N}{v}.
		\end{align*}
	\end{proof}
	
	Using the surface gradient notation
	\bes
	\GradGm{u} = N_m \times (\Grad{u} \times N_m),
	\ees
	and the fact that $\Curl{\Grad{u}} = \Curl{\Grad{v}} = 0$,
	we find that
	\begin{align*}
		a(\Grad{u},\Grad{v}) 
		& = - k_0^2 \InnerOmega{\bepsilon_0 \Grad{u}}{\Grad{v}} \\
		& \quad
		- \InnerGammau{(i \omega \mu_0) T_u[\GradGu{u}]}{\GradGu{v}}
		- \InnerGammaw{(i \omega \mu_0) T_w[\GradGw{u}]}{\GradGw{v}},
	\end{align*}
	and
	\bes
	L[\Grad{v}] = \InnerOmega{(\bepsilon_0 k_0^2 F)}{\Grad{v}}
	+ \InnerGammau{Q}{\GradGu{v}}
	+ \InnerGammaw{R}{\GradGw{v}}.
	\ees
	Using Lemma~\ref{Lemma:IBP} we discover
	\begin{align*}
		a(\Grad{u},\Grad{v}) 
		& = k_0^2 \InnerOmega{\Div{\bepsilon_0 \Grad{u}}}{v} \\
		& \quad
		- k_0^2 \InnerGammau{\bepsilon_0 \Grad{u} \cdot N_u}{v}
		- k_0^2 \InnerGammaw{\bepsilon_0 \Grad{u} \cdot N_w}{v} \\
		& \quad
		+ \InnerGammau{\DivGu{(i \omega \mu_0) T_u[\GradGu{u}]}}{v}
		+ \InnerGammaw{\DivGw{(i \omega \mu_0) T_w[\GradGw{u}]}}{v},
	\end{align*}
	and
	\begin{multline*}
		L[\Grad{v}] = -\InnerOmega{\Div{(\bepsilon_0 k_0^2 F)}}{v}
		+ \InnerGammau{(\bepsilon_0 k_0^2 F) \cdot N_u}{v}
		+ \InnerGammaw{(\bepsilon_0 k_0^2 F) \cdot N_w}{v}
		\\
		- \InnerGammau{\DivGu{Q}}{v}
		- \InnerGammaw{\DivGw{R}}{v}.
	\end{multline*}
	In this way (c.f.\ Lemma~3.24 of \cite{BaoLi22}) we see that
	$a(\Grad{u},\Grad{v}) = L(\Grad{v})$ is the weak
	formulation of the elliptic problem
	\begin{align*}
		& k_0^2 \Div{\bepsilon_0 \Grad{u}} = -\Div{(\bepsilon_0 k_0^2 F)}, 
		&& \Omega, \\
		& -k_0^2 \bepsilon_0 \Grad{u} \cdot N_u
		+ \DivGu{(i \omega \mu_0) T_u[\GradGu{u}]}
		= (\bepsilon_0 k_0^2 F) \cdot N - \DivGu{Q}, && \Gamma_u, \\
		& -k_0^2 \bepsilon_0 \Grad{u} \cdot N_w 
		+ \DivGw{(i \omega \mu_0) T_w[\GradGw{u}]}
		= (\bepsilon_0 k_0^2 F) \cdot N - \DivGw{R}, && \Gamma_w,
	\end{align*}
	which has a unique solution $u \in H_0^1(\Omega)$ \cite{BaoLi22}
	that satisfies
	\begin{multline*}
		\Norm{u}{H^1_0} 
		\leq C_e \left\{ \Norm{\Div{(\bepsilon_0 k_0^2 F)}}{L^2}
		+ \SobNorm{(\bepsilon_0 k_0^2 F) \cdot N_u}{-1/2} 
		+ \SobNorm{\DivGu{Q}}{-1/2}
		\right. \\ \left.
		+ \SobNorm{(\bepsilon_0 k_0^2 F) \cdot N_w}{-1/2}
		+ \SobNorm{\DivGw{R}}{-1/2}
		\right\},
	\end{multline*}
	or
	\be
	\label{Eqn:Hperp:Est}
	\Norm{u}{H^1_0}
	\leq C_e \left\{ \HdivNorm{(\bepsilon_0 k_0^2 F)}
	+ \HmhdivNorm{Q} + \HmhdivNorm{R} \right\}.
	\ee
	
	%
	%
	
	\subsection{Finding a Solution II: $\HH$}
	
	We continue by considering
	\bes
	a(\vec{u},\vec{v}) + a(\Grad{u},\vec{v}) = L[\vec{v}],
	\ees
	which, as we have recovered $u$, we rewrite as
	\bes
	a(\vec{u},\vec{v}) = L[\vec{v}] - a(\Grad{u},\vec{v}).
	\ees
	Bao \& Li \cite{BaoLi22} (Theorem~3.28) demonstrate
	that, save for a countable number of frequencies $\omega$
	(alternatively wavenumbers $k_0$),
	there exists a unique solution $\vec{u} \in \Hcurl$.
	Furthermore, from the inf--sup estimate of Bao \& Li
 \cite{BaoLi22}
	(Equation~(4.83))
	\bes
	\sup_{0 \neq \vec{v} \in \Hcurl}
	\frac{a(\vec{u},\vec{v})}{\HcurlNorm{\vec{v}}}
	\geq \gamma_1 \HcurlNorm{\vec{u}},
	\quad
	\forall\ \vec{u} \in \Hcurl,
	\ees
	for some $\gamma_1 > 0$, we have the inequality
	\bes
	a(\vec{u},\vec{v}) \geq \gamma_1 \HcurlNorm{\vec{u}}
	\HcurlNorm{\vec{v}}.
	\ees
	Therefore, from the continuity estimate established
	in the proof of Theorem~3.28 from Bao \& Li \cite{BaoLi22},
	\bes
	\Abs{a(\vec{u},\vec{v})}
	\leq C \HcurlNorm{\vec{u}} \HcurlNorm{\vec{v}},
	\ees
	we can show that
	\begin{align*}
		\gamma_1 \HcurlNorm{\vec{u}} \HcurlNorm{\vec{v}}
		& \leq \Abs{a(\vec{u},\vec{v})}
		\leq \Abs{L[\vec{v}]} + \Abs{a(\Grad{u},\vec{v})} \\
		& \leq \Abs{\InnerOmega{(\bepsilon_0 k_0^2 F)}{\vec{v}}}
		+ \Abs{\InnerGammau{Q}{\vec{v}}}
		+ \Abs{\InnerGammaw{R}{\vec{v}}} \\
		& \quad 
		+ C \HcurlNorm{\Grad{u}} \HcurlNorm{\vec{v}}.
	\end{align*}
	From Lemma~3.15 of \cite{BaoLi22} we have
	\bes
	\Abs{\InnerGammam{\vec{u}}{\vec{v}}}
	\leq C \HmhdivNorm{\vec{u}} \HmhdivNorm{\vec{v}},
	\ees
	so that
	\begin{align*}
		\gamma_1 \HcurlNorm{\vec{u}} \HcurlNorm{\vec{v}}
		& \leq \Norm{(\bepsilon_0 k_0^2 F)}{L^2} \Norm{\vec{v}}{L^2} \\
		& \quad
		+ \left( \HmhdivNorm{Q} + \HmhdivNorm{R} \right)
		\HmhcurlNorm{\vec{v}} \\
		& \quad
		+ C \Norm{u}{H^1_0} \HcurlNorm{\vec{v}}.
	\end{align*}
	Lemma~3.16 of \cite{BaoLi22} establishes that, for a given
	$\gamma_0 > 0$,
	\bes
	\HmhcurlNorm{\vec{u}}\leq \gamma_0 \HcurlNorm{\vec{u}},
	\ees
	so that
	\begin{multline*}
		\gamma_1 \HcurlNorm{\vec{u}} \HcurlNorm{\vec{v}} \\
		\leq \left\{ \Norm{(\bepsilon_0 k_0^2 F)}{L^2}
		+ \gamma_0 \left( \HmhdivNorm{Q}+ \HmhdivNorm{R} \right)
		+ C \Norm{u}{H^1_0} \right\} \HcurlNorm{\vec{v}}.
	\end{multline*}
	Clearly, by cancelling $\HcurlNorm{\vec{v}}$, there exists
	a constant $C_e > 0$ such that
	\bes
	\HcurlNorm{\vec{u}} \leq C_e
	\left\{ \HdivNorm{(\bepsilon_0 k_0^2 F)} 
	+ \HmhdivNorm{Q} + \HmhdivNorm{R} \right\}.
	\ees
	By combining this with \eqref{Eqn:Hperp:Est} we find
	\be
	\label{Eqn:Hcurl:Est}
	\HcurlNorm{E} \leq C_e
	\left\{ \HdivNorm{(\bepsilon_0 k_0^2 F)} 
	+ \HmhdivNorm{Q} + \HmhdivNorm{R} \right\}.
	\ee
	
	Next, simply applying the $L^2$ norm to both sides of 
	\eqref{Eqn:EllEst:b} we obtain
	\bes
	\Norm{\Div{\bepsilon_0 k_0^2 E}}{L^2} 
	= \Norm{\Div{\bepsilon k_0^2 F}}{L^2}.
	\ees
	Since we now know that $E \in \Hcurl$ we have $E \in L^2$, and,
	since $\bepsilon_0 \in L^{\infty}$,
	\bes
	\Norm{\bepsilon_0 k_0^2 E}{L^2} 
	\leq \SupNorm{\bepsilon_0} k_0^2 \Norm{E}{L^2}
	\leq \SupNorm{\bepsilon_0} k_0^2 \HcurlNorm{E}.
	\ees
	Now, recalling the definition of the $\Hdiv$ norm,
	\begin{align*}
		\HdivNorm{\bepsilon_0 k_0^2 E}^2
		& = \Norm{\bepsilon_0 k_0^2 E}{L^2}^2
		+ \Norm{\Div{\bepsilon_0 k_0^2 E}}{L^2}^2 \\
		& \leq \SupNorm{\bepsilon_0}^2 k_0^4 \HcurlNorm{E}^2
		+ \HdivNorm{\bepsilon_0 k_0^2 F}^2.
	\end{align*}
	Estimate \eqref{Eqn:Hcurl:Est} delivers \eqref{Eqn:EllEst:Est}.
	
	%
	%
	
	\section{The Proof of Theorem~\ref{Thm:JointAnal:EllEst}}
	\label{Sec:JointAnal:Proof}
	
	To establish Theorem~\ref{Thm:JointAnal:EllEst} we work by induction
	on the orders of spatial derivatives, beginning with the
	$x$--derivative, moving to the $y$--derivative, and concluding
	with the $z$--derivative. To begin this project we start with the
	following result on analyticity in $x$.
	\begin{theorem}
		\label{Thm:JointAnal:EllEst:x}
		Given any integer $q \geq 0$, if
		$(\omega,\bepsilon_0) \in \cP$,
		$\bepsilon_0 \in C^{\omega}_q(\Omega)$, such that
		\bes
		\HolderNorm{\frac{\px^r}{r!} \bepsilon_0}{q}
		\leq C_{\bepsilon_0} \frac{\eta^r}{(r+1)^2},
		\ees
		for all $r \geq 0$ and some constants 
		$C_{\bepsilon_0}, \eta > 0$,
		and $(\bepsilon_0 k_0^2 F) \in C^{\omega}(\Omega)$ such that
		\bes
		\HdivNorm{\frac{\px^r}{r!} \left[ \bepsilon_0 k_0^2 F \right]}
		\leq C_F \frac{\eta^r}{(r+1)^2},
		\ees
		for all $r \geq 0$ and some constant $C_F > 0$, and
		$Q \in C^{\omega}(\Gamma_u)$ and
		$R \in C^{\omega}(\Gamma_w)$ satisfying 
		\begin{align*}
			& \HmhdivNorm{\frac{\px^r}{r!} Q}
			\leq C_Q \frac{\eta^r}{(r+1)^2}, \\
			& \HmhdivNorm{\frac{\px^r}{r!} R}
			\leq C_R \frac{\eta^r}{(r+1)^2},
		\end{align*}
		for all $r \geq 0$ and some constants $C_R, C_Q > 0$.
		Then, there exists a unique solution
		$E \in C^{\omega}(\Omega)$ of 
		\begin{align*}
			& \cL_0 E = \bepsilon_0 k_0^2 F,
			&& \text{in $\Omega$}, \\
			& -\Div{\bepsilon_0 k_0^2 E} = \Div{ \bepsilon_0 k_0^2 F}, 
			&& \text{in $\Omega$}, \\
			& \cB_u E = Q, && \text{at $\Gamma_u$}, \\
			& \cB_w E = R, && \text{at $\Gamma_w$}, \\
			& E(x+d_x,y+d_y, z) = e^{i \alpha d_x + i \beta d_y} E(x,y,z),
		\end{align*}
		satisfying 
		\be
		\label{Eqn:JointAnal:EllEst:Estimate:x}
		\XNorm{ \frac{\px^r}{r!} E } 
		\leq \underline{C}_e \frac{\eta^r}{(r+1)^2},
		\ee
		for all $r \geq 0$ where
		\bes
		\underline{C}_e = \overline{C} (C_F + C_Q + C_R) > 0,
		\ees
		and $\overline{C} > 0$ is a constant.
	\end{theorem}
	\begin{proof}
		We prove this theorem by induction on $r \geq 0$. When $r = 0$
		we apply Theorem~\ref{Thm:EllEst} to conclude that
		\bes
		\XNorm{E} \leq C_e \left( \HdivNorm{\bepsilon_0 k_0^2 F}
		+ \HmhdivNorm{Q} + \HmhdivNorm{R} \right)
		\leq \underline{C}_e.
		\ees
		We now assume \eqref{Eqn:JointAnal:EllEst:Estimate:x} for
		all $r < \overline{r}$ and seek to establish this estimate
		at $r = \bar{r}$. Applying the differential operator
		$\px^r/r!$ to \eqref{Eqn:JointAnal:EllEst} delivers
		\begin{align*}
			& \cL_0 \left[ \frac{\px^r}{r!} E \right] 
			= \bepsilon_0 k_0^2 X_r,
			&& \text{in $\Omega$}, \\
			& -\Div{\bepsilon_0 k_0^2 \frac{\px^r}{r!} E} 
			= \Div{ \bepsilon_0 k_0^2 X_r},
			&& \text{in $\Omega$}, \\
			& \cB_u \left[ \frac{\px^r}{r!} E \right] 
			= \frac{\px^r}{r!} Q, && \text{at $\Gamma_u$}, \\
			& \cB_w \left[ \frac{\px^r}{r!} E \right] 
			= \frac{\px^r}{r!} R, && \text{at $\Gamma_w$}, \\
			& \frac{\px^r}{r!} E(x+d_x,y+d_y, z) 
			= e^{i \alpha d_x + i \beta d_y} \frac{\px^r}{r!} E(x,y,z),
		\end{align*}
		where
		\bes
		X_r := \frac{1}{\bepsilon_0 k_0^2} \left\{
		\frac{\px^r}{r!} \left[ \bepsilon_0 k_0^2 F \right]
		+ \left[ \cL_0, \frac{\px^r}{r!} \right] E
		\right\},
		\ees
		and $[A,B] := AB - BA$ is the commutator. With this we can
		express
		\bes
		A B = B A + [A,B] = B A - [B,A].
		\ees
		From Theorem~\ref{Thm:EllEst} we have
		\bes
		\XNorm{\frac{\px^{\bar{r}}}{\bar{r}!} E}
		\leq C_e \left\{
		\HdivNorm{\bepsilon_0 k_0^2 X_{\bar{r}}} 
		+ \HmhdivNorm{\frac{\px^{\bar{r}}}{\bar{r}!} Q}
		+ \HmhdivNorm{\frac{\px^{\bar{r}}}{\bar{r}!} R} \right\},
		\ees
		while Lemma~\ref{Lemma:XEst} gives
		\bes
		\XNorm{\frac{\px^{\bar{r}}}{\bar{r}!} E}
		\leq C_e \left\{
		C_F \frac{\eta^{\bar{r}}}{(\bar{r}+1)^2}
		+ \underline{C}_e \tilde{C} \frac{\eta^{\bar{r}-1}}{(\bar{r}+1)^2}
		+ C_Q \frac{\eta^{\bar{r}}}{(\bar{r}+1)^2}
		+ C_R \frac{\eta^{\bar{r}}}{(\bar{r}+1)^2} \right\}.
		\ees
		We are done provided that
		\bes
		\underline{C}_e \geq 2 C_e (C_F + C_Q + C_R),
		\quad
		\eta \geq 2 C_e \tilde{C}.
		\ees
	\end{proof}
	
	\begin{lemma}
		\label{Lemma:XEst}
		Given any integer $q \geq 0$, if
		$(\omega,\bepsilon_0) \in \cP$,
		$\bepsilon_0 \in C^{\omega}_q(\Omega)$, such that
		\bes
		\HolderNorm{\frac{\px^r}{r!} \bepsilon_0}{q}
		\leq C_{\bepsilon_0} \frac{\eta^r}{(r+1)^2},
		\ees
		for all $r \geq 0$ and some constants 
		$C_{\bepsilon_0}, \eta > 0$,
		and $(\bepsilon_0 k_0^2 F) \in C^{\omega}(\Omega)$ such that
		\bes
		\HdivNorm{\frac{\px^r}{r!} \left[ \bepsilon_0 k_0^2 F \right]}
		\leq C_F \frac{\eta^r}{(r+1)^2},
		\ees
		for all $r \geq 0$ and some constant $C_F > 0$. Assuming
		\bes
		\XNorm{ \frac{\px^r}{r!} E } 
		\leq \underline{C}_e \frac{\eta^r}{(r+1)^2},
		\ees
		for all $r < \bar{r}$ then
		\bes
		\HdivNorm{\bepsilon_0 k_0^2 X_{\bar{r}}}
		\leq 
		C_F \frac{\eta^{\bar{r}}}{(\bar{r}+1)^2}
		+ \underline{C}_e \tilde{C} \frac{\eta^{\bar{r}-1}}{(\bar{r}+1)^2}.
		\ees
	\end{lemma}
	\begin{proof}
		We note that
		\bes
		X_{\bar{r}} := \frac{1}{\bepsilon_0 k_0^2} \left\{
		\frac{\px^{\bar{r}}}{\bar{r}!} 
		\left[ \bepsilon_0 k_0^2 F \right]
		+ \left[ \cL_0, \frac{\px^{\bar{r}}}{\bar{r}!} \right] E
		\right\},
		\ees
		so that,
		\bes
		\HdivNorm{\bepsilon_0 k_0^2 X_{\bar{r}}}
		\leq \HdivNorm{\frac{\px^{\bar{r}}}{\bar{r}!}
			\left[ \bepsilon_0 k_0^2 F \right]}
		+ \HdivNorm{\left[ \cL_0, \frac{\px^{\bar{r}}}{\bar{r}!} \right] E}.
		\ees
		The first term can be bounded above by
		\bes
		C_F \frac{\eta^{\bar{r}}}{(\bar{r}+1)^2},
		\ees
		by assumption. Regarding the second term, we have, for any $r \geq 0$,
		\begin{align*}
			\left[ \cL_0, \frac{\px^r}{r!} \right] E
			& = \cL_0 \left[ \frac{\px^r}{r!} E \right] 
			- \frac{\px^r}{r!} \left[ \cL_0 [E] \right] \\
			& = \Curl{\Curl{\frac{\px^r}{r!} E}} 
			- \bepsilon_0 k_0^2 \frac{\px^r}{r!} E
			- \frac{\px^r}{r!} \left[ \Curl{\Curl{E}} 
			- \bepsilon_0 k_0^2 E \right] \\
			& = - \bepsilon_0 k_0^2 \frac{\px^r}{r!} E
			+ \frac{\px^r}{r!} \left[ \bepsilon_0 k_0^2 E \right].
		\end{align*}
		The Leibniz rule tells us that
		\begin{align*}
			\left[ \cL_0, \frac{\px^r}{r!} \right] E
			& = - \bepsilon_0 k_0^2 \frac{\px^r}{r!} E
			+ k_0^2 \sum_{j=0}^{r} 
			\left( \frac{\px^{r-j}}{(r-j)!} \bepsilon_0 \right)
			\frac{\px^j}{j!} E \\
			& = k_0^2 \sum_{j=0}^{r-1} 
			\left( \frac{\px^{r-j}}{(r-j)!} \bepsilon_0 \right)
			\frac{\px^j}{j!} E.
		\end{align*}
		Setting $r = \bar{r}$ we can estimate
		\begin{align*}
			\HdivNorm{\left[ \cL_0, \frac{\px^{\bar{r}}}{\bar{r}!} \right] E}
			& \leq k_0^2 \sum_{j=0}^{\bar{r}-1} 
			\HdivNorm{ \left( \frac{\px^{\bar{r}-j}}{(\bar{r}-j)!} \bepsilon_0 \right)
				\frac{\px^j}{j!} E } \\
			& \leq k_0^2 M \sum_{j=0}^{\bar{r}-1} 
			\HolderNorm{ \frac{\px^{\bar{r}-j}}{(\bar{r}-j)!} \bepsilon_0 }{1}
			\HdivNorm{ \frac{\px^j}{j!} E } \\
			& \leq k_0^2 M \sum_{j=0}^{\bar{r}-1} 
			\frac{1}{\bar{r}-j}
			\HolderNorm{ \frac{\px^{\bar{r}-j-1}}{(\bar{r}-j-1)!} \bepsilon_0 }{2}
			\HdivNorm{ \frac{\px^j}{j!} \left[ E \right] } \\
			& \leq k_0^2 M \sum_{j=0}^{\bar{r}-1}
			C_{\bepsilon_0} \frac{\eta^{\bar{r}-j-1}}{(\bar{r}-j-1+1)^2}
			\underline{C}_e \frac{\eta^j}{(j+1)^2} \\
			& \leq k_0^2 M C_{\bepsilon_0} \underline{C}_e S
			\frac{\eta^{\bar{r}-1}}{(\bar{r}+1)^2},
		\end{align*}
		and we choose
		\bes
		\tilde{C} \geq k_0^2 M C_{\bepsilon_0} S.
		\ees
	\end{proof}
	
	The next step is to prove analyticity jointly in $x$ and $y$.
	\begin{theorem}
		\label{Thm:JointAnal:EllEst:xy}
		Given any integer $q \geq 0$, if
		$(\omega,\bepsilon_0) \in \cP$,
		$\bepsilon_0 \in C^{\omega}_q(\Omega)$, such that
		\bes
		\HolderNorm{\frac{\px^r \py^t}{(r+t)!} \bepsilon_0}{q}
		\leq C_{\bepsilon_0} \frac{\eta^r}{(r+1)^2}
		\frac{\theta^t}{(t+1)^2},
		\ees
		for all $r, t \geq 0$ and some constants 
		$C_{\bepsilon_0}, \eta, \theta > 0$,
		and $(\bepsilon_0 k_0^2 F) \in C^{\omega}(\Omega)$ such that
		\bes
		\HdivNorm{\frac{\px^r \py^t}{(r+t)!} \left[ \bepsilon_0 k_0^2 F \right]}
		\leq C_F \frac{\eta^r}{(r+1)^2} \frac{\theta^t}{(t+1)^2},
		\ees
		for all $r, t \geq 0$ and some constant $C_F > 0$, and
		$Q \in C^{\omega}(\Gamma_u)$ and
		$R \in C^{\omega}(\Gamma_w)$ satisfying 
		\begin{align*}
			& \HmhdivNorm{\frac{\px^r \py^t}{(r+t)!} Q}
			\leq C_Q \frac{\eta^r}{(r+1)^2} \frac{\theta^t}{(t+1)^2}, \\
			& \HmhdivNorm{\frac{\px^r \py^t}{(r+t)!} R}
			\leq C_R \frac{\eta^r}{(r+1)^2} \frac{\theta^t}{(t+1)^2},
		\end{align*}
		for all $r, t \geq 0$ and some constants $C_R, C_Q > 0$.
		Then, there exists a unique solution
		$E \in C^{\omega}(\Omega)$ of 
		\begin{align*}
			& \cL_0 E = \bepsilon_0 k_0^2 F,
			&& \text{in $\Omega$}, \\
			& -\Div{\bepsilon_0 k_0^2 E} = \Div{ \bepsilon_0 k_0^2 F}, 
			&& \text{in $\Omega$}, \\
			& \cB_u E = Q, && \text{at $\Gamma_u$}, \\
			& \cB_w E = R, && \text{at $\Gamma_w$}, \\
			& E(x+d_x,y+d_y, z) = e^{i \alpha d_x + i \beta d_y} E(x,y,z),
		\end{align*}
		satisfying 
		\be
		\label{Eqn:JointAnal:EllEst:Estimate:xy}
		\XNorm{ \frac{\px^r \py^t}{(r+t)!} E } 
		\leq \underline{C}_e \frac{\eta^r}{(r+1)^2}
		\frac{\theta^t}{(t+1)^2},
		\ee
		for all $r,t \geq 0$ where
		\bes
		\underline{C}_e = \overline{C} (C_F + C_Q + C_R) > 0,
		\ees
		and $\overline{C} > 0$ is a constant.
	\end{theorem}
	\begin{proof}
		We prove this theorem by induction on $t \geq 0$. When $t = 0$
		we apply Theorem~\ref{Thm:JointAnal:EllEst:x} to establish
		\eqref{Eqn:JointAnal:EllEst:Estimate:xy}.
		We now assume \eqref{Eqn:JointAnal:EllEst:Estimate:xy} for
		all $r \geq 0$ and $t < \overline{t}$ and seek to establish this estimate
		at $t = \bar{t}$. Applying the differential operator
		$(\px^r \py^t)/(r+t)!$ to \eqref{Eqn:JointAnal:EllEst} delivers
		\begin{align*}
			& \cL_0 \left[ \frac{\px^r \py^t}{(r+t)!} E \right] 
			= \bepsilon_0 k_0^2 Y_{r,t},
			&& \text{in $\Omega$}, \\
			& -\Div{\bepsilon_0 k_0^2 \frac{\px^r \py^t}{(r+t)!} E} 
			= \Div{ \bepsilon_0 k_0^2 Y_{r,t}},
			&& \text{in $\Omega$}, \\
			& \cB_u \left[ \frac{\px^r \py^t}{(r+t)!} E \right] 
			= \frac{\px^r \py^t}{(r+t)!} Q, && \text{at $\Gamma_u$}, \\
			& \cB_w \left[ \frac{\px^r \py^t}{(r+t)!} E \right] 
			= \frac{\px^r \py^t}{(r+t)!} R, && \text{at $\Gamma_w$}, \\
			& \frac{\px^r \py^t}{(r+t)!} E(x+d_x,y+d_y, z) 
			= e^{i \alpha d_x + i \beta d_y} \frac{\px^r \py^t}{(r+t)!} E(x,y,z),
		\end{align*}
		where
		\bes
		Y_{r,t} := \frac{1}{\bepsilon_0 k_0^2} \left\{
		\frac{\px^r \py^t}{(r+t)!} \left[ \bepsilon_0 k_0^2 F \right]
		+ \left[ \cL_0, \frac{\px^r \py^t}{(r+t)!} \right] E
		\right\}.
		\ees
		From Theorem~\ref{Thm:EllEst} we have
		\begin{multline*}
			\XNorm{\frac{\px^r \py^{\bar{t}}}{(r+\bar{t})!} E}
			\leq C_e \left\{
			\HdivNorm{\bepsilon_0 k_0^2 Y_{r,\bar{t}}} 
			\right. \\ \left.
			+ \HmhdivNorm{\frac{\px^r \py^{\bar{t}}}{(r+\bar{t})!} Q}
			+ \HmhdivNorm{\frac{\px^r \py^{\bar{t}}}{(r+\bar{t})!} R} \right\},
		\end{multline*}
		while Lemma~\ref{Lemma:XYEst} gives
		\begin{multline*}
			\XNorm{\frac{\px^r \px^{\bar{t}}}{\bar{r}!} E}
			\leq C_e \left\{
			C_F \frac{\eta^r}{(r+1)^2} \frac{\theta^{\bar{t}}}{(\bar{t}+1)^2}
			+ \underline{C}_e \tilde{C} \frac{\eta^r}{(r+1)^2}
			\frac{\theta^{\bar{t}-1}}{(\bar{t}+1)^2}
			\right. \\ \left.
			+ C_Q \frac{\eta^r}{(r+1)^2} \frac{\theta^{\bar{t}}}{(\bar{t}+1)^2}
			+ C_R \frac{\eta^r}{(r+1)^2} \frac{\theta^{\bar{t}}}{(\bar{t}+1)^2}
			\right\}.
		\end{multline*}
		We are done provided that
		\bes
		\underline{C}_e \geq 2 C_e (C_F + C_Q + C_R),
		\quad
		\theta \geq 2 C_e \tilde{C}.
		\ees
	\end{proof}
	
	\begin{lemma}
		\label{Lemma:XYEst}
		Given any integer $q \geq 0$, if
		$(\omega,\bepsilon_0) \in \cP$,
		$\bepsilon_0 \in C^{\omega}_q(\Omega)$, such that
		\bes
		\HolderNorm{\frac{\px^r \py^t}{(r+t)!} \bepsilon_0}{q}
		\leq C_{\bepsilon_0} \frac{\eta^r}{(r+1)^2}
		\frac{\theta^t}{(t+1)^2},
		\ees
		for all $r, t \geq 0$ and some constants 
		$C_{\bepsilon_0}, \eta, \theta > 0$,
		and $(\bepsilon_0 k_0^2 F) \in C^{\omega}(\Omega)$ such that
		\bes
		\HdivNorm{\frac{\px^r \py^t}{(r+t)!} \left[ \bepsilon_0 k_0^2 F \right]}
		\leq C_F \frac{\eta^r}{(r+1)^2} \frac{\theta^t}{(t+1)^2},
		\ees
		for all $r, t \geq 0$ and some constant $C_F > 0$. Assuming
		\bes
		\XNorm{ \frac{\px^r \py^t}{(r+t)!} E } 
		\leq \underline{C}_e \frac{\eta^r}{(r+1)^2}
		\frac{\theta^t}{(t+1)^2},
		\ees
		for all $r \geq 0$ and $t < \bar{t}$ then
		\bes
		\HdivNorm{\bepsilon_0 k_0^2 Y_{r,\bar{t}}}
		\leq 
		C_F \frac{\eta^r}{(r+1)^2} \frac{\theta^{\bar{t}}}{(\bar{t}+1)^2}
		+ \underline{C}_e \tilde{C} \frac{\eta^r}{(r+1)^2} 
		\frac{\theta^{\bar{t}-1}}{(\bar{t}+1)^2}.
		\ees
	\end{lemma}
	\begin{proof}
		Very similar to that of Lemma~\ref{Lemma:XEst} and therefore
		omitted.
	\end{proof}
	
We are now in a position to prove 
Theorem~\ref{Thm:JointAnal:EllEst}. Once again, we work
by induction, this time on $s$, the order of the $z$
derivative acting on $E$. As before, at $s=0$ we use
Theorem~\ref{Thm:JointAnal:EllEst:xy} to establish
\eqref{Eqn:JointAnal:EllEst:Estimate}.
We now assume \eqref{Eqn:JointAnal:EllEst:Estimate},
\be
\label{inductive_hypothesis}
\XNorm{ \frac{\px^r \py^t \pz^s}{(r+t+s)!} E } 
  \leq \underline{C}_e \frac{\eta^r}{(r+1)^2}
  \frac{\theta^t}{(t+1)^2} \frac{\zeta^s}{(s+1)^2},
\ee
for all $r, t \geq 0$ and all $s < \bar{s}$. We now
examine this estimate at order $s = \bar{s}$. 
From the definition of the $X$--norm we have
\begin{align*}
\XNorm{ \frac{\px^r \py^t \pz^{\bar{s}}}{(r+t+\bar{s})!} E }^2
  & = \HcurlNorm{ \frac{\px^r \py^t \pz^{\bar{s}}}{(r+t+\bar{s})!} E }^2
  + \HdivNorm{ \bepsilon_0 k_0^2 \frac{\px^r \py^t \pz^{\bar{s}}}{(r+t+\bar{s})!} E }^2
\\
  & = \Norm{ \frac{\px^r \py^t \pz^{\bar{s}}}{(r+t+\bar{s})!} E }{L^2}^2
  + \Norm{ \Curl{ \frac{\px^r \py^t \pz^{\bar{s}}}{(r+t+\bar{s})!} E } }{L^2}^2 \\
  & \quad 
  + \Norm{ \bepsilon_0 k_0^2
  \frac{\px^r \py^t \pz^{\bar{s}}}{(r+t+\bar{s})!} E }{L^2}^2
  + \Norm{ \Div{ \bepsilon_0 k_0^2 
  \frac{\px^r \py^t \pz^{\bar{s}}}{(r+t+\bar{s})!} E } }{L^2}^2 \\
  & = \Norm{ \frac{\px^r \py^t \pz^{\bar{s}}}{(r+t+\bar{s})!} E }{L^2}^2
  + \Norm{ \frac{\px^r \py^t \pz^{\bar{s}}}{(r+t+\bar{s})!} \Curl{ E } }{L^2}^2 \\
  & \quad 
  + \Norm{ \bepsilon_0 k_0^2
  \frac{\px^r \py^t \pz^{\bar{s}}}{(r+t+\bar{s})!} E }{L^2}^2 \\
  & \quad
  + \Norm{ \Grad{\bepsilon_0} \cdot \left( k_0^2
  \frac{\px^r \py^t \pz^{\bar{s}}}{(r+t+\bar{s})!} E \right) }{L^2}^2
  + \Norm{ \bepsilon_0 k_0^2
  \frac{\px^r \py^t \pz^{\bar{s}}}{(r+t+\bar{s})!} \Div{ E } }{L^2}^2 \\
  & \leq \left\{ 1 + \HolderNorm{\bepsilon_0}{1}^2 k_0^4 \right\} 
  \Norm{ \frac{\px^r \py^t \pz^{\bar{s}}}{(r+t+\bar{s})!} E }{L^2}^2 \\
  & \quad
  + \Norm{ \frac{\px^r \py^t \pz^{\bar{s}}}{(r+t+\bar{s})!} \Curl{ E } }{L^2}^2
  + \Norm{\bepsilon_0 k_0^2 
  \frac{\px^r \py^t \pz^{\bar{s}}}{(r+t+\bar{s})!} \Div{ E } }{L^2}^2 \\
  & =: \left\{ 1 + \HolderNorm{\bepsilon_0}{1}^2 k_0^4 \right\} I_1 
  + I_2 + I_3.
\end{align*}

\textbf{The estimate of $I_1$:}
Using the notation $E = (E^x,E^y,E^z)^T$ we obtain
\begin{align*}
I_1 & = \Norm{ \frac{\px^r \py^t \pz^{\bar{s}}}{(r+t+\bar{s})!} E^x }{L^2}^2
  + \Norm{ \frac{\px^r \py^t \pz^{\bar{s}}}{(r+t+\bar{s})!} E^y }{L^2}^2
  + \Norm{ \frac{\px^r \py^t \pz^{\bar{s}}}{(r+t+\bar{s})!} E^z }{L^2}^2 \\
  & = \Norm{ \frac{\px^r \py^t \pz^{\bar{s}-1}}{(r+t+\bar{s})!} \pz E^x }{L^2}^2
  + \Norm{ \frac{\px^r \py^t \pz^{\bar{s}-1}}{(r+t+\bar{s})!} \pz E^y }{L^2}^2
  +  \Norm{ \frac{\px^r \py^t \pz^{\bar{s}-1}}{(r+t+\bar{s})!} \pz E^z }{L^2}^2 \\
  &	:= Z_1 + Z_2 + Z_3.
\end{align*}
To estimate $Z_1$ and $Z_2$, we notice that 
\begin{align*}
\Curl{E} & = \left( \py E^z - \pz E^y, \pz E^x - \px E^z,
  \px E^y - \py E^x \right) \\
  & = \left( \py E^z, -\px E^z, \px E^y - \py E^x \right)
  + \left( -\pz E^y, \pz E^x , 0 \right),
\end{align*}
we obtain 
\bes
\left( -\pz E^y, \pz E^x , 0 \right)
  = \Curl{E} - \left( \py E^z , -\px E^z, \px E^y - \py E^x \right).
\ees
Therefore,
\begin{align*}
\left\{ Z_1 + Z_2 \right\}^{1/2} 
  & =  \Norm{ \frac{\px^r \py^t \pz^{\bar{s}-1}}{(r+t+\bar{s})!}
  (-\pz E^y, \pz E^x, 0)^T}{L^2} \\
  & \leq \Norm{\frac{\px^r \py^t \pz^{\bar{s}-1}}
  {(r+t+\bar{s})!} \Curl{E}}{L^2} \\
  & \quad +
  \Norm{\frac{\px^r \py^t \pz^{\bar{s}-1}}{(r+t+\bar{s})!}
  (\py E^z, -\px E^z, \px E^y - \py E^x)}{L^2} \\
  & \leq \Norm{\frac{\px^r \py^t \pz^{\bar{s}-1}}{(r+t+\bar{s})!} \Curl{E}}{L^2} \\
  & \quad 
  + \Norm{\frac{\px^r \py^{t+1} \pz^{\bar{s}-1}}{(r+t+\bar{s})!} E^z}{L^2} 
  + \Norm{\frac{\px^{r+1} \py^{t} \pz^{\bar{s}-1}}{(r+t+\bar{s})!} E^z}{L^2} \\
  & \quad 
  + \Norm{\frac{\px^{r+1} \py^{t} \pz^{\bar{s}-1}}{(r+t+\bar{s})!} E^y}{L^2} 
  + \Norm{\frac{\px^r \py^{t+1} \pz^{\bar{s}-1}}{(r+t+\bar{s})!} E^x}{L^2}.
\end{align*}
Using the inequality $(r+t+\bar{s}-1)! \leq (r+t+\bar{s})!$
and the inductive hypothesis, \eqref{inductive_hypothesis},
we have
\begin{align*}
\left\{ Z_1 + Z_2 \right\}^{1/2} 
  & \leq \underline{C}_e \frac{\eta^r}{(r+1)^2}
  \frac{\theta^t}{(t+1)^2} \frac{\zeta^{\bar{s}-1}}{(\bar{s}-1+1)^2} \\
  & \quad 
  + 2 \underline{C}_e \frac{\eta^r}{(r+1)^2}
  \frac{\theta^{t+1}}{(t+1+1)^2} \frac{\zeta^{\bar{s}-1}}{(\bar{s}-1+1)^2} \\
  & \quad
  + 2 \underline{C}_e \frac{\eta^{r+1}}{(r+1+1)^2}
  \frac{\theta^{t}}{(t+1)^2} \frac{\zeta^{\bar{s}-1}}{(\bar{s}-1+1)^2} \\
  & = \underline{C}_e \left\{ \frac{1}{\zeta} \frac{(\bar{s}+1)^2}{\bar{s}^2} 
  + 2 \theta \frac{(t+1)^2}{(t+2)^2} 
  \frac{1}{\zeta} \frac{(\bar{s}+1)^2}{\bar{s}^2} 
  + 2 \eta \frac{(r+1)^2}{(r+2)^2} 
  \frac{1}{\zeta} \frac{(\bar{s}+1)^2}{\bar{s}^2} \right\} \\
  & \quad 
  \times \frac{\eta^r}{(r+1)^2} \frac{\theta^t}{(t+1)^2}
  \frac{\zeta^{\bar{s}}}{(\bar{s}+1)^2}.
\end{align*}
With this and using the inequalities
\begin{align*}
  & \frac{(\bar{s}+1)^2}{\bar{s}^2} \leq 4, && \bar{s} \geq 1, \\
  & \frac{(t+1)^2}{(t+2)^2} < 1, && t \geq 0,
\end{align*}
we obtain
\begin{align*}
\left\{ Z_1 + Z_2 \right\}^{1/2} 
  & \leq \underline{C}_e \left\{ \frac{4 + 8 \theta + 8 \eta}{\zeta} \right\}
  \frac{\eta^r}{(r+1)^2} \frac{\theta^t}{(t+1)^2} 
  \frac{\zeta^{\bar{s}}}{(\bar{s}+1)^2}, \\
  & \leq \underline{C}_e \frac{\eta^r}{(r+1)^2} \frac{\theta^t}{(t+1)^2}
  \frac{\zeta^{\bar{s}}}{(\bar{s}+1)^2},
\end{align*}
where $\zeta \geq (4 + 8\theta + 8\eta)$.
Next, we estimate $Z_3$ by using the fact that 
\bes
\pz E^z = \Div{E} - \px E^x - \py E^y.
\ees
Therefore, 
\begin{align*}
Z_3^{1/2} 
  & = \Norm{ \frac{\px^r \py^t \pz^{\bar{s}-1}}{(r+t+\bar{s})!} 
  (0,0,\pz E^z)^T }{L^2} \\
  & \leq \Norm{\frac{\px^r \py^t \pz^{\bar{s}-1}}{(r+t+\bar{s})!} \Div{E}}{L^2}
  + \Norm{\frac{\px^{r+1} \py^t \pz^{\bar{s}-1}}{(r+t+\bar{s})!}E^x}{L^2}
  + \Norm{\frac{\px^r \py^{t+1} \pz^{\bar{s}-1}}{(r+t+\bar{s})!}E^y}{L^2} \\
  & =: W_1 + W_2 + W_3.
\end{align*}
While $W_2$ and $W_3$ can be directly controlled by using the inductive 
hypothesis \eqref{inductive_hypothesis}, $W_1$ can be estimated as follows:
\begin{align*}
W_1 & = \Norm{ \frac{\bepsilon_0 k^2_0}{\bepsilon_0 k_0^2}
  \frac{\px^r\py^t\pz^{\bar{s}-1}}{(r+t+\bar{s})!}\Div{E}}{L^2} \\
  & \leq \frac{1}{\min \Abs{\bepsilon_0} k_0^2}
  \Norm{\bepsilon_0 k_0^2 \frac{\px^r\py^t\pz^{\bar{s}-1}}{(r+t+\bar{s})!}
  \Div{E}}{L^2} \\
  & \leq \frac{1}{\min \Abs{\bepsilon_0} k_0^2}
  \frac{(r+t+\bar{s}-1)!}{(r+t+\bar{s})!}
  \frac{(\bar{s}+1)^2}{\bar{s}^2 \zeta}
  \underline{C}_e \frac{\eta^r}{(r+1)^2}
  \frac{\theta^t}{(t+1)^2} \frac{\zeta^{\bar{s}}}{(\bar{s}+1)^2} \\
  & \leq \frac{1}{\min \Abs{\bepsilon_0} k_0^2}
  \frac{4}{\zeta} \underline{C}_e \frac{\eta^r}{(r+1)^2}
  \frac{\theta^t}{(t+1)^2} \frac{\zeta^{\bar{s}}}{(\bar{s}+1)^2} \\
  & \leq \underline{C}_e \frac{\eta^r}{(r+1)^2}
  \frac{\theta^t}{(t+1)^2} \frac{\zeta^{\bar{s}}}{(\bar{s}+1)^2},
\end{align*}
where $\zeta > 4/(\min \Abs{\bepsilon_0} k_0^2)$, which completes our
estimation of $I_1$.

\textbf{The estimate of $I_2$:} We begin with the fact that
\bes
\Curl{\pz E} = \left( \py \pz E^z - \pz^2 E^y, \pz^2 E^x - \px\pz E^z,
  \px\pz E^y - \py\pz E^x \right).
\ees
So,
\begin{align*}
I_2^{1/2} & = \Norm{\frac{\px^r \py^t \pz^{\overline{s}-1}}
  {(r+t+\overline{s})!} \pz \Curl{E}}{L^2} \\
  & = \Norm{\frac{\px^r \py^t \pz^{\overline{s}-1}}{(r+t+\overline{s})!}
  \Curl{\pz E}}{L^2} \\
  & \leq \Norm{\frac{\px^r \py^{t+1} \pz^{\overline{s}-1}}{(r+t+\overline{s})!}
  \pz E^z}{L^2} 
  + \Norm{\frac{\px^{r+1} \py^t \pz^{\overline{s}-1}}
  {(r+t+\overline{s})!}\pz E^z}{L^2} \\
  & \quad 
  + \Norm{\frac{\px^{r+1} \py^{t} \pz^{\overline{s}-1}}
  {(r+t+\overline{s})!} \pz E^y}{L^2} 
  + \Norm{\frac{\px^{r} \py^{t+1} \pz^{\overline{s}-1}}
  {(r+t+\overline{s})!} \pz E^x}{L^2} \\
  & \quad 
  + \Norm{\frac{\px^{r} \py^{t} \pz^{\overline{s}-1}}
  {(r+t+\overline{s})!} \pz^2 E^y}{L^2} 
  + \Norm{\frac{\px^{r} \py^{t} \pz^{\overline{s}-1}}{(r+t+\overline{s})!}
  \pz^2 E^x}{L^2} \\
  & =: X_1 + X_2 + X_3 + X_4 + X_5 + X_6.
\end{align*}
The terms $X_1$ and $X_2$ are similar to $Z_3$ while the terms $X_3$ and $X_4$
can be controlled by using the same techniques as shown in bounding
$Z_1$ and $Z_2$. Since $X_5$ and $X_6$ are quite similar to estimate, we fix on
the term $X_5$. 
		
To begin this estimation, we can write the governing equation,
$\mathcal{L}_0 E = \bepsilon_0 k_0^2F$, as
\bes
-\Laplacian{E} + \Grad{\Div{E}} - k_0^2 \bepsilon_0 E = \bepsilon_0 k_0^2F, 
\ees
which implies that
\bes
\pz^2 E^y = -\bepsilon_0 k_0^2 F^y - k_0^2 \bepsilon_0 E^y
  - \px^2 E^y + \px \py E^x + \py \pz E^z.
\ees
Therefore, 
\begin{align*}
X_5 & = \Norm{\frac{\px^{r} \py^{t} \pz^{\overline{s}-1}}
  {(r+t+\overline{s})!} \pz^2 E^y}{L^2} \\
  & \leq \Norm{\frac{\px^{r} \py^{t} \pz^{\overline{s}-1}}
  {(r+t+\overline{s})!} (\bepsilon_0 k_0^2 F^y) }{L^2} 
  + \Norm{\frac{\px^{r} \py^{t} \pz^{\overline{s}-1}}{(r+t+\overline{s})!}
  (\bepsilon_0 k_0^2 E^y)}{L^2} 
  + \Norm{\frac{\px^{r+2} \py^{t} \pz^{\overline{s}-1}}{(r+t+\overline{s})!} E^y}{L^2} \\
  & \quad 
  + \Norm{\frac{\px^{r+1} \py^{t+1} \pz^{\overline{s}-1}}
  {(r+t+\overline{s})!} E^x}{L^2} 
  + \Norm{\frac{\px^{r} \py^{t+1} \pz^{\overline{s}-1}}{(r+t+\overline{s})!}
  \pz E^z}{L^2} \\
  & := P_1 + P_2 + P_3 + P_4 + P_5.
\end{align*}
It is clear that $P_1$ can be estimated by using the analyticity assumption
on $\bepsilon_0 k_0^2 F$, and the terms $P_3$ and $P_4$ can be controlled
by the inductive hypothesis \eqref{inductive_hypothesis}. On the other hand,
the term $P_5$ can be estimated by using the same techniques we applied to
$Z_3$. It remains to bound $P_2$. Using Leibniz's rule, for any $r ,t, s \geq 0$,
we have
\begin{align*}
\frac{\px^r \py^t \pz^{s}}{(r+t+s)!}(\bepsilon_0 k_0^2 E) 
  & = \frac{r! t! s!}{(r+t+s)!} \sum_{j=0}^{r} \sum_{k=0}^{t} \sum_{\ell=0}^{s}
  \left( \frac{\px^{r-j}}{(r-j)!} \frac{\py^{t-k}}{(t-k)!}
  \frac{\pz^{s-\ell}}{(s-\ell)^2} \bepsilon_0 k_0^2 \right)
\\
  & \quad 
  \times \left( \frac{\px^j}{j!} \frac{\py^k}{k!} \frac{\pz^{\ell}}{\ell !} E \right).
\end{align*}
Using the triangle inequality, the fact that 
$r! t! (\overline{s}-1)! \leq (r+t+\overline{s})!$,
the hypotheses on $\bepsilon_0$, and the inductive hypothesis
\eqref{inductive_hypothesis} we obtain
\begin{align*}
P_2 & = \Norm{\frac{\px^r \py^t \pz^{\overline{s}-1}}{(r+t+\overline{s})!}
  (\bepsilon_0k_0^2E)}{L^2} \\
  & \leq \frac{r! t! (\overline{s}-1)!}{(r+t+\overline{s})!}
  \sum_{j=0}^{r} \sum_{k=0}^{t}\sum_{\ell =0}^{\overline{s}-1}
  \Norm{ \left( \frac{\px^{r-j}}{(r-j)!} \frac{\py^{t-k}}{(t-k)!} 
  \frac{\pz^{\overline{s}-1-\ell}}{(\overline{s}-1-\ell)^2} 
  \bepsilon_0 k_0^2 \right)
  \left( \frac{\px^j}{j!} \frac{\py^k}{k!} \frac{\pz^{\ell}}{\ell !} E \right)}{L^2} \\
  & \leq \sum_{j=0}^{r} \sum_{k=0}^{t} \sum_{\ell =0}^{{\overline{s}-1}}
  \SupNorm{ \frac{\px^{r-j}}{(r-j)!} \frac{\py^{t-k}}{(t-k)!}
  \frac{\pz^{{\overline{s}-1}-\ell}}{(\overline{s}-1-\ell)!}
  \bepsilon_0k_0^2 }
  \Norm{\frac{\px^j}{j!} \frac{\py^k}{k!} \frac{\pz^{\ell}}{\ell!} E}{L^2} \\
  & \leq \sum_{j=0}^{r} \sum_{k=0}^{t} \sum_{\ell=0}^{{\overline{s}-1}}
  C_{\bepsilon_0} \frac{\eta^{r-j}}{(r-j+1)^2}
  \frac{\theta^{t-k}}{(t-k+1)^2}
  \frac{\zeta^{\overline{s}-1-\ell}}{(\overline{s}-1-\ell +1)^2} \\
  & \quad \times \underline{C}_e 
  \frac{\eta^j}{(j+1)^2} \frac{\theta^k}{(k+1)^2}\frac{\zeta^{\ell}}{(\ell+1)^2} \\
  & \leq \underline{C}_e C_{\bepsilon_0} \frac{\eta^r}{(r+1)^2}
  \frac{\theta^t}{(t+1)^2} \frac{\zeta^{\overline{s}-1}}{(\overline{s}-1+1)^2} \\
  & \quad \times \sum_{j=0}^{r}\sum_{k=0}^{t} \sum_{\ell =0}^{{\overline{s}-1}} 
  \frac{(r+1)^2(t+1)^2}{(r-j+1)^2(j+1)^2(t-k+1)^2(k+1)^2}
  \frac{(\overline{s}-1 +1)^2}{(\overline{s}-1-\ell + 1)^2(\ell +1)^2} \\
  & \leq \underline{C}_e C_{\bepsilon_0} S^3
  \frac{(\overline{s}+1)^2}{\overline{s}^2\zeta} \frac{\eta^r}{(r+1)^2}
  \frac{\theta^t}{(t+1)^2} \frac{\zeta^{\overline{s}}}{(\overline{s}+1)^2} \\
  & \leq \underline{C}_e \frac{\eta^r}{(r+1)^2}\frac{\theta^t}{(t+1)^2}
  \frac{\zeta^{\overline{s}}}{(\overline{s}+1)^2},
\end{align*}
for some $\zeta \geq 4 C_{\bepsilon_0}S^3$. 

\textbf{The estimate of $I_3$:} We conclude with the computation
\begin{align*}
I_3^{1/2} 
  & = \Norm{\bepsilon_0 k_0^2 \frac{\px^r \py^t \pz^{\bar{s}}}{(r+t+\bar{s})!}
  \Div{ E } }{L^2} \\
  & \leq \Norm{\bepsilon_0 k_0^2 \frac{\px^r \py^t \pz^{\bar{s}}}
  {(r+t+\bar{s})!} \px E^x }{L^2} 
  + \Norm{\bepsilon_0 k_0^2 \frac{\px^r \py^t \pz^{\bar{s}}}
  {(r+t+\bar{s})!} \py E^y }{L^2} \\
  & \quad
  + \Norm{\bepsilon_0 k_0^2 \frac{\px^r \py^t \pz^{\bar{s}}}
  {(r+t+\bar{s})!} \pz E^z }{L^2} \\
  & \leq \Norm{\bepsilon_0 k_0^2 \frac{\px^{r+1} \py^t \pz^{\bar{s}-1}}
  {(r+t+\bar{s})!} \pz E^x }{L^2}
  + \Norm{\bepsilon_0 k_0^2 \frac{\px^r \py^{t+1} \pz^{\bar{s}-1}}
  {(r+t+\bar{s})!} \pz E^y }{L^2} \\
  & \quad
  + \Norm{\bepsilon_0 k_0^2 \frac{\px^r \py^t \pz^{\bar{s}}}
  {(r+t+\bar{s})!} \pz E^z }{L^2} \\
  & =: L_1 + L_2 + L_3.
\end{align*}
The terms $L_1$ and $L_2$ are similar to $Z_1$ and $Z_2$, and we estimate
them in the same fashion. We estimate $L_3$ in the following way. We begin
by noting that $\Div{\bepsilon_0 k_0^2 E} = -\Div{ \bepsilon_0 k_0^2 F}$
implies that
\bes
\pz (\bepsilon_0 k_0^2 E^z) = -\Div{\bepsilon_0 k_0^2 F} 
  - \px(\bepsilon_0 k_0^2 E^x) - \py(\bepsilon_0k_0^2 E^y),
\ees
which also gives
\be
\label{B4}
\bepsilon_0 k_0^2 \pz E^z = -\Div{\bepsilon_0 k_0^2 F}
  - \px(\bepsilon_0 k_0^2 E^x) - \py(\bepsilon_0 k_0^2 E^y)
  - \pz(\bepsilon_0k_0^2) E^z.
\ee
In anticipation of our future estimation of higher derivatives
of this term we use the commutator notation to express
\bes
\frac{\px^r \py^t \pz^{\bar{s}}}{(r+t+\bar{s})!} 
  \left[ (\bepsilon_0 k_0^2) \pz E^z \right] 
  = \bepsilon_0 k_0^2 \frac{\px^r \py^t \pz^{\bar{s}}}{(r+t+\bar{s})!}
  \left[ \pz E^z \right] 
  + \left[ \frac{\px^r \py^t \pz^{\bar{s}}}{(r+t+\bar{s})!},
  \bepsilon_0 k_0^2 \right] \pz E^z,
\ees
where, by Leibniz's Rule,
\begin{align*}
\left[ \frac{\px^r \py^t \pz^{\bar{s}}}{(r+t+\bar{s})!},
  \bepsilon_0 k_0^2 \right] \pz E^z
  & = \frac{r! t! \bar{s}!}{(r+t+\bar{s}!)}
  \sum_{j=0}^{r} \sum_{k=0}^{t} \sum_{\ell=0}^{\bar{s}}
  \left( \frac{\px^{r-j}}{(r-j)!} \frac{\py^{t-k}}{(t-k)!}
  \frac{\pz^{\bar{s}-\ell}}{(\bar{s}-\ell)!} 
  \left[ \bepsilon_0 k_0^2 \right] \right)
  \\ & \quad \times
  \left( \frac{\px^j}{j!} \frac{\py^k}{k!}
  \frac{\pz^{\ell}}{\ell!} \left[ \pz E^z \right] \right)
  - (\bepsilon_0 k_0^2) \frac{\px^r \py^t \pz^{\bar{s}}}
  {(r+t+\bar{s})!} \pz E^z
  \\
  & = \frac{r! t! \bar{s}!}{(r+t+\bar{s}!)}
  \sum_{j=0}^{r} \sum_{k=0}^{t} \sum_{\ell=0}^{\bar{s}-1}
  \left( \frac{\px^{r-j}}{(r-j)!} \frac{\py^{t-k}}{(t-k)!}
  \frac{\pz^{\bar{s}-\ell}}{(\bar{s}-\ell)!} 
  \left[ \bepsilon_0 k_0^2 \right] \right)
  \\ & \quad \times
  \left( \frac{\px^j}{j!} \frac{\py^k}{k!}
  \frac{\pz^{\ell}}{\ell!} \left[ \pz E^z \right] \right)
  \\ & \quad
  + \frac{r! t! \bar{s}!}{(r+t+\bar{s}!)}
  \sum_{j=0}^{r} \sum_{k=0}^{t}
  \left( \frac{\px^{r-j}}{(r-j)!} \frac{\py^{t-k}}{(t-k)!} 
  \left[ \bepsilon_0 k_0^2 \right] \right)
  \\ & \quad \times
  \left( \frac{\px^j}{j!} \frac{\py^k}{k!}
  \frac{\pz^{\bar{s}}}{\bar{s}!} \left[ \pz E^z \right] \right)
  - (\bepsilon_0 k_0^2) \frac{\px^r \py^t \pz^{\bar{s}}}
  {(r+t+\bar{s})!} \pz E^z.
\end{align*}
Continuing
\begin{align*}
\left[ \frac{\px^r \py^t \pz^{\bar{s}}}{(r+t+\bar{s})!},
  \bepsilon_0 k_0^2 \right] \pz E^z
  & = \frac{r! t! \bar{s}!}{(r+t+\bar{s}!)}
  \sum_{j=0}^{r} \sum_{k=0}^{t} \sum_{\ell=0}^{\bar{s}-1}
  \left( \frac{\px^{r-j}}{(r-j)!} \frac{\py^{t-k}}{(t-k)!}
  \frac{\pz^{\bar{s}-\ell}}{(\bar{s}-\ell)!} 
  \left[ \bepsilon_0 k_0^2 \right] \right)
  \\ & \quad \times
  \left( \frac{\px^j}{j!} \frac{\py^k}{k!}
  \frac{\pz^{\ell}}{\ell!} \left[ \pz E^z \right] \right)
  \\ & \quad
  + \frac{r! t! \bar{s}!}{(r+t+\bar{s}!)}
  \sum_{j=0}^{r} \sum_{k=0}^{t-1}
  \left( \frac{\px^{r-j}}{(r-j)!} \frac{\py^{t-k}}{(t-k)!} 
  \left[ \bepsilon_0 k_0^2 \right] \right)
  \\ & \quad \times
  \left( \frac{\px^j}{j!} \frac{\py^k}{k!}
  \frac{\pz^{\bar{s}}}{\bar{s}!} \left[ \pz E^z \right] \right)
  \\ & \quad
  + \frac{r! t! \bar{s}!}{(r+t+\bar{s}!)}
  \sum_{j=0}^{r}
  \left( \frac{\px^{r-j}}{(r-j)!} 
  \left[ \bepsilon_0 k_0^2 \right] \right)
  \\ & \quad \times
  \left( \frac{\px^j}{j!} \frac{\py^t}{t!}
  \frac{\pz^{\bar{s}}}{\bar{s}!} \left[ \pz E^z \right] \right)
  - (\bepsilon_0 k_0^2) \frac{\px^r \py^t \pz^{\bar{s}}}
  {(r+t+\bar{s})!} \pz E^z,
\end{align*}
and finally,
\begin{align*}
\left[ \frac{\px^r \py^t \pz^{\bar{s}}}{(r+t+\bar{s})!},
  \bepsilon_0 k_0^2 \right] \pz E^z
  & = \frac{r! t! \bar{s}!}{(r+t+\bar{s}!)}
  \sum_{j=0}^{r} \sum_{k=0}^{t} \sum_{\ell=0}^{\bar{s}-1}
  \left( \frac{\px^{r-j}}{(r-j)!} \frac{\py^{t-k}}{(t-k)!}
  \frac{\pz^{\bar{s}-\ell}}{(\bar{s}-\ell)!} 
  \left[ \bepsilon_0 k_0^2 \right] \right)
  \\ & \quad \times
  \left( \frac{\px^j}{j!} \frac{\py^k}{k!}
  \frac{\pz^{\ell}}{\ell!} \left[ \pz E^z \right] \right)
  \\ & \quad
  + \frac{r! t! \bar{s}!}{(r+t+\bar{s}!)}
  \sum_{j=0}^{r} \sum_{k=0}^{t-1}
  \left( \frac{\px^{r-j}}{(r-j)!} \frac{\py^{t-k}}{(t-k)!} 
  \left[ \bepsilon_0 k_0^2 \right] \right)
  \\ & \quad \times
  \left( \frac{\px^j}{j!} \frac{\py^k}{k!}
  \frac{\pz^{\bar{s}}}{\bar{s}!} \left[ \pz E^z \right] \right)
  \\ & \quad
  + \frac{r! t! \bar{s}!}{(r+t+\bar{s}!)}
  \sum_{j=0}^{r-1}
  \left( \frac{\px^{r-j}}{(r-j)!} 
  \left[ \bepsilon_0 k_0^2 \right] \right)
  \\ & \quad \times
  \left( \frac{\px^j}{j!} \frac{\py^t}{t!}
  \frac{\pz^{\bar{s}}}{\bar{s}!} \left[ \pz E^z \right] \right).
\end{align*}
With this commutator notation, \eqref{B4} implies that
\begin{align*}
\bepsilon_0 k_0^2 \frac{\px^r\py^t\pz^{\bar{s}}}{(r+t+\bar{s})!} \pz E^z
  & = \frac{\px^r \py^t \pz^{\bar{s}}}{(r+t+\bar{s})!} 
  \left[ (\bepsilon_0 k_0^2) \pz E^z \right] 
  - \left[ \frac{\px^r \py^t \pz^{\bar{s}}}{(r+t+\bar{s})!},
  \bepsilon_0 k_0^2 \right] \pz E^z \\
  & = -\frac{\px^r \py^t \pz^{\bar{s}}}{(r+t+\bar{s})!}
  \Div{\bepsilon_0k_0^2 F} 
  - \frac{\px^{r+1} \py^t \pz^{\bar{s}}}{(r+t+\bar{s})!}
  \left[ \bepsilon_0k_0^2 E^x \right] \\
  & \quad
  - \frac{\px^r \py^{t+1} \pz^{\bar{s}}}{(r+t+\bar{s})!}
  \left[ \bepsilon_0 k_0^2 E^y \right] 
  - \frac{\px^r \py^t \pz^{\bar{s}}}{(r+t+\bar{s})!}
  \left[ \pz(\bepsilon_0 k_0^2) E^z \right] \\
  & \quad
  - \left[ \frac{\px^r \py^t \pz^{\bar{s}}}{(r+t+\bar{s})!},
  \bepsilon_0 k_0^2 \right] \pz E^z.
\end{align*}
Proceeding, we find
\begin{align*}
\bepsilon_0 k_0^2 \frac{\px^r\py^t\pz^{\bar{s}}}{(r+t+\bar{s})!} \pz E^z
  & = -\frac{\px^r \py^t \pz^{\bar{s}}}{(r+t+\bar{s})!} 
  \Div{\bepsilon_0 k_0^2 F} 
  - \frac{\px^{r+1} \py^t \pz^{\bar{s}-1}}{(r+t+\bar{s})!} 
  \pz \left[ \bepsilon_0 k_0^2 E^x \right] \\
  & \quad
  - \frac{\px^r \py^{t+1} \pz^{\bar{s}-1}}{(r+t+\bar{s})!}
  \pz \left[ \bepsilon_0 k_0^2 E^y \right] 
  - \frac{\px^r \py^t \pz^{\bar{s}-1}}{(r+t+\bar{s})!}
  \pz \left[ \pz(\bepsilon_0 k_0^2) E^z \right] \\
  & \quad
  - \left[ \frac{\px^r \py^t \pz^{\bar{s}}}{(r+t+\bar{s})!},
  \bepsilon_0 k_0^2 \right] \pz E^z.
\end{align*}
Using the product rule we find
\begin{align*}
\bepsilon_0 k_0^2 \frac{\px^r\py^t\pz^{\bar{s}}}{(r+t+\bar{s})!} \pz E^z
  & = - \frac{\px^r \py^t\pz^{\bar{s}}}{(r+t+\bar{s})!}
  \Div{\bepsilon_0 k_0^2 F} 
  - \frac{\px^{r+1} \py^t \pz^{\bar{s}-1}}{(r+t+\bar{s})!}
  \left[ \pz(\bepsilon_0 k_0^2) E^x \right] \\
  & \quad
  -  \frac{\px^{r+1} \py^t \pz^{\bar{s}-1}}{(r+t+\bar{s})!}
  \left[ \bepsilon_0 k_0^2 \pz E^x \right] 
  - \frac{\px^r \py^{t+1} \pz^{\bar{s}-1}}{(r+t+\bar{s})!}
  \left[ \pz(\bepsilon_0 k_0^2) E^y \right] \\
  & \quad
  - \frac{\px^r \py^{t+1} \pz^{\bar{s}-1}}{(r+t+\bar{s})!}
  \left[ \bepsilon_0 k_0^2 \pz E^y \right] 
  - \frac{\px^r \py^t \pz^{\bar{s}-1}}{(r+t+\bar{s})!}
  \left[ \pz^2(\bepsilon_0k_0^2) E^z \right] \\
  & \quad 
  - \frac{\px^r \py^t \pz^{\bar{s}-1}}{(r+t+\bar{s})!}
  \left[ \pz(\bepsilon_0 k_0^2 ) \pz E^z \right] \\
  & \quad
  - \left[ \frac{\px^r \py^t \pz^{\bar{s}}}{(r+t+\bar{s})!},
  \bepsilon_0 k_0^2 \right] \pz E^z.
\end{align*}
With this, we can estimate $L_3$ as follows:
\begin{align*}
L_3 & = \Norm{\bepsilon_0 k_0^2 \frac{\px^r \py^t \pz^{\bar{s}}}
  {(r+t+\bar{s})!} \pz E^z }{L^2} \\
  & \leq \Norm{\frac{\px^r \py^t \pz^{\bar{s}}}{(r+t+\bar{s})!} \Div{\bepsilon_0 k_0^2 F} }{L^2} \\
  & \quad +
  \Norm{ \frac{\px^{r+1} \py^t \pz^{\bar{s}-1}}
  {(r+t+\bar{s})!} \left[ \pz(\bepsilon_0 k_0^2) E^x \right]}{L^2}
  + \Norm{ \frac{\px^{r+1} \py^t \pz^{\bar{s}-1}}{(r+t+\bar{s})!}
  \left[ \bepsilon_0 k_0^2 \pz E^x \right] }{L^2} \\
  & \quad +
  \Norm{\frac{\px^r \py^{t+1} \pz^{\bar{s}-1}}{(r+t+\bar{s})!}
  \left[ \pz(\bepsilon_0 k_0^2) E^y \right]}{L^2}
  + \Norm{\frac{\px^r \py^{t+1} \pz^{\bar{s}-1}}{(r+t+\bar{s})!}
  \left[ (\bepsilon_0 k_0^2) \pz E^y \right]}{L^2} \\
  & \quad +
  \Norm{\frac{\px^r \py^t \pz^{\bar{s}-1}}{(r+t+\bar{s})!}
  \left[ \pz^2(\bepsilon_0 k_0^2) E^z \right]}{L^2} 
  + \Norm{\frac{\px^r \py^t \pz^{\bar{s}-1}}{(r+t+\bar{s})!}
  \left[ \pz(\bepsilon_0 k_0^2) \pz E^z \right]}{L^2} \\
  & \quad + 
  \Norm{\left[ \frac{\px^r \py^t \pz^{\bar{s}}}{(r+t+\bar{s})!},
  \bepsilon_0 k_0^2 \right] \pz E^z}{L^2} \\
  & =: U_1 + U_2 + U_3 + U_4 + U_5 + U_6 + U_7 + U_8.
\end{align*}
It is clear that $U_1$ can be bounded by using the analyticity 
assumption on $\bepsilon_0 k_0^2 F$. Next, the terms $U_2, U_4, U_6$ 
can be readily controlled by using the same techniques we used in
bounding $P_2$, with the help of the Leibniz's rule and the inductive
hypothesis \eqref{inductive_hypothesis}. On the other hand, 
$U_3$ and $U_5$ are similar, so we simply present the estimation of
$U_3$. (We turn to $U_7$ and $U_8$ in a moment.)

Using the Leibniz's rule, we obtain
\begin{multline*}
U_3 \leq \frac{r!t!\bar{s}!}{(r+t+\bar{s})!} 
  \sum_{j=0}^{r} \sum_{k=0}^{t} \sum_{\ell=0}^{\bar{s}-1}
  \SupNorm{\frac{\px^{r-j}}{(r-j)!} \frac{\py^{t-k}}{(t-k)!}
  \frac{\pz^{\bar{s}-1-\ell}}{(\bar{s}-1-\ell)^2} 
  \left[ \bepsilon_0 k_0^2 \right]}
  \\ \times
  \Norm{\frac{\px^j}{j!} \frac{\py^k}{k!} \frac{\pz^{\ell}}{\ell!}
  \pz E^x}{L^2}.
\end{multline*}
Following the approach of $Z_1$, we obtain 
\begin{align*}
\Norm{\frac{\px^j}{j!} \frac{\py^k}{k!} \frac{\pz^{\ell}}{\ell!}
  \pz E^x}{L^2}
  & \leq \Norm{\frac{\px^j}{j!} \frac{\py^k}{k!}
  \frac{\pz^{\ell}}{\ell!} \Curl{E}}{L^2} 
  + \Norm{\frac{\px^j}{j!} \frac{\py^{k+1}}{k!} 
  \frac{\pz^{\ell}}{\ell!} E^z}{L^2} \\
  & \quad
  + \Norm{\frac{\px^{j+1}}{j!} \frac{\py^k}{k!}
  \frac{\pz^{\ell}}{\ell!} E^z}{L^2} 
  + \Norm{\frac{\px^{j+1}}{j!} \frac{\py^k}{k!}
  \frac{\pz^{\ell}}{\ell!} E^y}{L^2} \\
  & \quad
  + \Norm{\frac{\px^j}{j!} \frac{\py^{k+1}}{k!}
  \frac{\pz^{\ell}}{\ell!} E^x}{L^2}.
\end{align*}
With this, we can bound $U_3$ as follows
\begin{multline*}
U_3 \leq \sum_{j=0}^{r} \sum_{k=0}^{t} \sum_{\ell=0}^{\bar{s}-1}
  \SupNorm{\frac{\px^{r-j}}{(r-j)!} \frac{\py^{t-k}}{(t-k)!}
  \frac{\pz^{\bar{s}-1-\ell}}{(\bar{s}-1-\ell)^2} 
  \left[ \bepsilon_0 k_0^2 \right]} \\
  \times \left\{ \Norm{\frac{\px^j}{j!} \frac{\py^k}{k!} 
  \frac{\pz^{\ell}}{\ell!} \Curl{E}}{L^2} 
  + \Norm{\frac{\px^j}{j!} \frac{\py^{k+1}}{k!} 
  \frac{\pz^{\ell}}{\ell!} E^z}{L^2} 
  + \Norm{\frac{\px^{j+1}}{j!} \frac{\py^k}{k!}
  \frac{\pz^{\ell}}{\ell!} E^z}{L^2} 
  \right. \\ \left.
  + \Norm{\frac{\px^{j+1}}{j!} \frac{\py^k}{k!}
  \frac{\pz^{\ell}}{\ell!} E^y}{L^2} 
  + \Norm{\frac{\px^j}{j!} \frac{\py^{k+1}}{k!}
  \frac{\pz^{\ell}}{\ell!} E^x}{L^2}
  \right\}.
\end{multline*}
Using the analyticity assumption on $\bepsilon_0$ and the inductive
hypothesis \eqref{inductive_hypothesis}, we get
\begin{align*}
U_3	& \leq \sum_{j=0}^{r} \sum_{k=0}^{t} \sum_{\ell=0}^{{\overline{s}-1}}
  C_{\bepsilon_0} \frac{\eta^{r-j}}{(r-j+1)^2}
  \frac{\theta^{t-k}}{(t-k+1)^2}
  \frac{\zeta^{\overline{s}-1-\ell}}{(\overline{s}-1-\ell +1)^2}
  \\
  & \quad \times
  \underline{C}_e \left\{ \frac{\eta^j}{(j+1)^2}
  \frac{\theta^k}{(k+1)^2} \frac{\zeta^{\ell}}{(\ell+1)^2}
  + 2 \frac{\eta^j}{(j+1)^2} \frac{\theta^{k+1}}{(k+2)^2}
  \frac{\zeta^{\ell}}{(\ell+1)^2} \right. \\
  & \quad
  \left. + 2 \frac{\eta^{j+1}}{(j+2)^2} 
  \frac{\theta^k}{(k+1)^2} \frac{\zeta^{\ell}}{(\ell+1)^2} \right\} \\
  & \leq \sum_{j=0}^{r} \sum_{k=0}^{t} \sum_{\ell=0}^{{\overline{s}-1}}
  C_{\bepsilon_0} \frac{\eta^{r-j}}{(r-j+1)^2}
  \frac{\theta^{t-k}}{(t-k+1)^2}
  \frac{\zeta^{\overline{s}-1-\ell}}{(\overline{s}-1-\ell +1)^2}
  \\
  & \quad
  \times \underline{C}_e (1 + 2 \theta + 2 \eta) 
  \frac{\eta^j}{(j+1)^2} \frac{\theta^k}{(k+1)^2}
  \frac{\zeta^{\ell}}{(\ell+1)^2} \\
  & \leq \underline{C}_e C_{\bepsilon_0} (1 + 2 \theta + 2 \eta)
  \frac{\eta^r}{(r+1)^2} \frac{\theta^t}{(t+1)^2}
  \frac{\zeta^{\overline{s}-1}}{(\overline{s}-1+1)^2}
  \\
  & \quad \times \sum_{j=0}^{r} \sum_{k=0}^{t} 
  \sum_{\ell=0}^{{\overline{s}-1}} 
  \frac{(r+1)^2(t+1)^2}{(r-j+1)^2(j+1)^2(t-k+1)^2(k+1)^2}
  \frac{(\overline{s}-1 +1)^2}{(\overline{s}-1-\ell + 1)^2(\ell +1)^2}
  \\
  & \leq 
  \underline{C}_e C_{\bepsilon_0} (1 + 2 \theta + 2 \eta) 
  \frac{1}{\zeta} S^3 \frac{(\overline{s}+1)^2}{\overline{s}^2} 
  \frac{\eta^r}{(r+1)^2} \frac{\theta^t}{(t+1)^2} 
  \frac{\zeta^{\overline{s}}}{(\overline{s}+1)^2}
  \\
  & \leq \underline{C}_e \frac{\eta^r}{(r+1)^2}
  \frac{\theta^t}{(t+1)^2} 
  \frac{\zeta^{\overline{s}}}{(\overline{s}+1)^2},
\end{align*}
for some $\zeta \geq 4 C_{\bepsilon_0} (1 + 2 \theta + 2 \eta) S^3$.

Finally, we observe that, by using Leibniz's rule, $U_7$ and $U_8$
can be controlled in a similar fashion. So, it suffices to estimate
$U_8$, and for this we compute
\begin{align*}
U_8 & = 
  \Norm{ \left[ \frac{\px^r \py^t \pz^{\bar{s}}}{(r+t+\bar{s})!},
  \bepsilon_0 k_0^2 \right] \pz E^z }{L^2} \\
  & \leq \frac{r! t! \bar{s}!}{(r+t+\bar{s}!)}
  \sum_{j=0}^{r} \sum_{k=0}^{t} \sum_{\ell=0}^{\bar{s}-1}
  \SupNorm{ \frac{\px^{r-j}}{(r-j)!} \frac{\py^{t-k}}{(t-k)!}
  \frac{\pz^{\bar{s}-\ell}}{(\bar{s}-\ell)!} 
  \left[ \bepsilon_0 k_0^2 \right] }
  \Norm{ \frac{\px^j}{j!} \frac{\py^k}{k!}
  \frac{\pz^{\ell}}{\ell!} \left[ \pz E^z \right] }{L^2}
  \\ & \quad
  + \frac{r! t! \bar{s}!}{(r+t+\bar{s}!)}
  \sum_{j=0}^{r} \sum_{k=0}^{t-1}
  \SupNorm{ \frac{\px^{r-j}}{(r-j)!} \frac{\py^{t-k}}{(t-k)!} 
  \left[ \bepsilon_0 k_0^2 \right] }
  \Norm{ \frac{\px^j}{j!} \frac{\py^k}{k!}
  \frac{\pz^{\bar{s}}}{\bar{s}!} \left[ \pz E^z \right] }{L^2}
  \\ & \quad
  + \frac{r! t! \bar{s}!}{(r+t+\bar{s}!)}
  \sum_{j=0}^{r-1}
  \SupNorm{ \frac{\px^{r-j}}{(r-j)!} 
  \left[ \bepsilon_0 k_0^2 \right] }
  \Norm{ \frac{\px^j}{j!} \frac{\py^t}{t!}
  \frac{\pz^{\bar{s}}}{\bar{s}!} \left[ \pz E^z \right] }{L^2}.
\end{align*}
Using the inductive hypothesis, \eqref{inductive_hypothesis}, and
following the same approach as used to address $U_3$, we get
\bes
U_8 \leq \underline{C}_e \frac{\eta^r}{(r+1)^2}
  \frac{\theta^t}{(t+1)^2} \frac{\zeta^{\overline{s}}}{(\overline{s}+1)^2},
\ees
for some $\zeta \geq 4 C_{\bepsilon_0} (1 + \theta + \eta) S^3$, and
the proof is complete.

%
%
	
\bibliographystyle{abbrv}
\bibliography{nicholls}

\begin{thebibliography}{10}

\bibitem{BaoLi22}
G.~Bao and P.~Li.
\newblock {\em Maxwell's equations in periodic structures}, volume 208 of {\em
  Applied Mathematical Sciences}.
\newblock Springer, Singapore; Science Press Beijing, Beijing, [2022]
  \copyright 2022.

\bibitem{ColtonKress13}
D.~Colton and R.~Kress.
\newblock {\em Inverse acoustic and electromagnetic scattering theory},
  volume~93 of {\em Applied Mathematical Sciences}.
\newblock Springer, New York, third edition, 2013.

\bibitem{DevilleFischerMund02}
M.~O. Deville, P.~F. Fischer, and E.~H. Mund.
\newblock {\em High-order methods for incompressible fluid flow}, volume~9 of
  {\em Cambridge Monographs on Applied and Computational Mathematics}.
\newblock Cambridge University Press, Cambridge, 2002.

\bibitem{ELGTW98}
T.~W. Ebbesen, H.~J. Lezec, H.~F. Ghaemi, T.~Thio, and P.~A. Wolff.
\newblock Extraordinary optical transmission through sub-wavelength hole
  arrays.
\newblock {\em Nature}, 391(6668):667--669, 1998.

\bibitem{ESHES06}
I.~El-Sayed, X.~Huang, and M.~El-Sayed.
\newblock Selective laser photo-thermal therapy of epithelial carcinoma using
  anti-egfr antibody conjugated gold nanoparticles.
\newblock {\em Cancer Lett.}, 239(1):129--135, 2006.

\bibitem{ErnstGander12}
O.~G. Ernst and M.~J. Gander.
\newblock Why it is difficult to solve {H}elmholtz problems with classical
  iterative methods.
\newblock In {\em Numerical analysis of multiscale problems}, volume~83 of {\em
  Lect. Notes Comput. Sci. Eng.}, pages 325--363. Springer, Heidelberg, 2012.

\bibitem{Evans10}
L.~C. Evans.
\newblock {\em Partial differential equations}.
\newblock American Mathematical Society, Providence, RI, second edition, 2010.

\bibitem{FengLinLorton15}
X.~Feng, J.~Lin, and C.~Lorton.
\newblock An efficient numerical method for acoustic wave scattering in random
  media.
\newblock {\em SIAM/ASA J. Uncertain. Quantif.}, 3(1):790--822, 2015.

\bibitem{FengLinLorton16}
X.~Feng, J.~Lin, and C.~Lorton.
\newblock A multimodes {M}onte {C}arlo finite element method for elliptic
  partial differential equations with random coefficients.
\newblock {\em Int. J. Uncertain. Quantif.}, 6(5):429--443, 2016.

\bibitem{GottliebOrszag77}
D.~Gottlieb and S.~A. Orszag.
\newblock {\em Numerical analysis of spectral methods: theory and
  applications}.
\newblock Society for Industrial and Applied Mathematics, Philadelphia, Pa.,
  1977.
\newblock {C}BMS-NSF Regional Conference Series in Applied Mathematics, No. 26.

\bibitem{HesthavenWarburton08}
J.~S. Hesthaven and T.~Warburton.
\newblock {\em Nodal discontinuous {G}alerkin methods}, volume~54 of {\em Texts
  in Applied Mathematics}.
\newblock Springer, New York, 2008.
\newblock Algorithms, analysis, and applications.

\bibitem{Homola08}
J.~Homola.
\newblock Surface plasmon resonance sensors for detection of chemical and
  biological species.
\newblock {\em Chemical Reviews}, 108(2):462--493, 2008.

\bibitem{Ihlenburg98}
F.~Ihlenburg.
\newblock {\em Finite element analysis of acoustic scattering}.
\newblock Springer-Verlag, New York, 1998.

\bibitem{Johnson87}
C.~Johnson.
\newblock {\em Numerical solution of partial differential equations by the
  finite element method}.
\newblock Cambridge University Press, Cambridge, 1987.

\bibitem{JJJLWO13}
J.~Jose, L.~R. Jordan, T.~W. Johnson, S.~H. Lee, N.~J. Wittenberg, and S.-H.
  Oh.
\newblock Topographically flat substrates with embedded nanoplasmonic devices
  for biosensing.
\newblock {\em Adv Funct Mater}, 23:2812--2820, 2013.

\bibitem{Kress14}
R.~Kress.
\newblock {\em Linear integral equations}.
\newblock Springer-Verlag, New York, third edition, 2014.

\bibitem{LeVeque07}
R.~J. LeVeque.
\newblock {\em Finite difference methods for ordinary and partial differential
  equations}.
\newblock Society for Industrial and Applied Mathematics (SIAM), Philadelphia,
  PA, 2007.
\newblock Steady-state and time-dependent problems.

\bibitem{LJJOO12}
N.~C. Lindquist, T.~W. Johnson, J.~Jose, L.~M. Otto, and S.-H. Oh.
\newblock Ultrasmooth metallic films with buried nanostructures for backside
  reflection-mode plasmonic biosensing.
\newblock {\em Annalen der Physik}, 524:687--696, 2012.

\bibitem{MoiolaSpence14}
A.~Moiola and E.~A. Spence.
\newblock Is the {H}elmholtz equation really sign-indefinite?
\newblock {\em SIAM Rev.}, 56(2):274--312, 2014.

\bibitem{Moskovits85}
M.~Moskovits.
\newblock Surface--enhanced spectroscopy.
\newblock {\em Reviews of Modern Physics}, 57(3):783--826, 1985.

\bibitem{Nicholls19b}
D.~P. Nicholls.
\newblock A high--order perturbation of envelopes ({HOPE}) method for
  scattering by periodic inhomogeneous media.
\newblock {\em Quarterly of Applied Mathematics}, 78:725--757, 2020.

\bibitem{NichollsReitich99}
D.~P. Nicholls and F.~Reitich.
\newblock A new approach to analyticity of {D}irichlet-{N}eumann operators.
\newblock {\em Proc. Roy. Soc. Edinburgh Sect. A}, 131(6):1411--1433, 2001.

\bibitem{NichollsReitich00b}
D.~P. Nicholls and F.~Reitich.
\newblock Analytic continuation of {D}irichlet-{N}eumann operators.
\newblock {\em Numer. Math.}, 94(1):107--146, 2003.

\bibitem{NichollsReitich03a}
D.~P. Nicholls and F.~Reitich.
\newblock Shape deformations in rough surface scattering: Cancellations,
  conditioning, and convergence.
\newblock {\em J. Opt. Soc. Am. A}, 21(4):590--605, 2004.

\bibitem{NichollsReitich03b}
D.~P. Nicholls and F.~Reitich.
\newblock Shape deformations in rough surface scattering: Improved algorithms.
\newblock {\em J. Opt. Soc. Am. A}, 21(4):606--621, 2004.

\bibitem{NichollsReitichJohnsonOh14}
D.~P. Nicholls, F.~Reitich, T.~W. Johnson, and S.-H. Oh.
\newblock Fast high--order perturbation of surfaces ({HOPS}) methods for
  simulation of multi--layer plasmonic devices and metamaterials.
\newblock {\em Journal of the Optical Society of America, A}, 31(8):1820--1831,
  2014.

\bibitem{NichollsVo23}
D.~P. Nicholls and L.~Vo.
\newblock A high-order perturbation of envelopes ({HOPE}) method for vector
  electromagnetic scattering by periodic inhomogeneous media.
\newblock {\em (submitted)}, 2023.

\bibitem{Petit80}
R.~Petit, editor.
\newblock {\em Electromagnetic theory of gratings}.
\newblock Springer-Verlag, Berlin, 1980.

\bibitem{SauterSchwab11}
S.~A. Sauter and C.~Schwab.
\newblock {\em Boundary element methods}, volume~39 of {\em Springer Series in
  Computational Mathematics}.
\newblock Springer-Verlag, Berlin, 2011.
\newblock Translated and expanded from the 2004 German original.

\bibitem{ShenTang06}
J.~Shen and T.~Tang.
\newblock {\em Spectral and high-order methods with applications}, volume~3 of
  {\em Mathematics Monograph Series}.
\newblock Science Press Beijing, Beijing, 2006.

\bibitem{ShenTangWang11}
J.~Shen, T.~Tang, and L.-L. Wang.
\newblock {\em Spectral methods}, volume~41 of {\em Springer Series in
  Computational Mathematics}.
\newblock Springer, Heidelberg, 2011.
\newblock Algorithms, analysis and applications.

\bibitem{Strikwerda04}
J.~C. Strikwerda.
\newblock {\em Finite difference schemes and partial differential equations}.
\newblock Society for Industrial and Applied Mathematics (SIAM), Philadelphia,
  PA, second edition, 2004.

\bibitem{Yeh05}
P.~Yeh.
\newblock {\em Optical waves in layered media}, volume~61.
\newblock Wiley-Interscience, 2005.

\end{thebibliography}
	
\end{document}